\documentclass[12pt]{article}
\usepackage{amsmath}
\newcommand{\bee}{\begin{eqnarray*}}
\newcommand{\ene}{\end{eqnarray*}}
\newcommand{\beeq}{\begin{equation}}
\newcommand{\eneq}{\end{equation}}
\newtheorem{lem}{Lemma}[section]
\newcommand{\bel}{\begin{lem}}
\newcommand{\enl}{\end{lem}}
\newtheorem{defi}{Definition}[section]
\newcommand{\bef}{\begin{defi}}
\newcommand{\enf}{\end{defi}}
\newtheorem{exap}{Example}[section]
\newcommand{\beex}{\begin{exap}}
\newcommand{\enex}{\end{exap}}
\newtheorem{theo}{Theorem}[section]
\newcommand{\beth}{\begin{theo}}
\newcommand{\enth}{\end{theo}}
\newtheorem{prop}{Proposition}[section]
\newcommand{\bep}{\begin{prop}}
\newcommand{\enp}{\end{prop}}
\newtheorem{cor}{Corollary}[section]
\newcommand{\bec}{\begin{cor}}
\newcommand{\enc}{\end{cor}}
\newtheorem{rem}{Remark}[section]
\newcommand{\ber}{\begin{rem}}
\newcommand{\enr}{\end{rem}}

\usepackage{amsmath}
\usepackage{graphicx}
\usepackage{latexsym}
\usepackage{latexsym}
\usepackage{graphicx}

\begin{document}

\title{THE DERIVATIVE OF INFLUENCE FUNCTION,\\ LOCATION  BREAKDOWN POINT,\\
%GROUP  LEVERAGE \\
GROUP  INFLUENCE \\
AND\\
REGRESSION RESIDUALS' PLOTS
%\\
%IN REGRESSION  
% PLOTS OF RESIDUALS\\
%AND\\
% LEVERAGE 
%CASES
}
 \author {Yannis G. Yatracos\\
Faculty of Communication and Media Studies\\
Cyprus University of Technology}

%The National University of Singapore}
%\date{December 9, 2004}
\maketitle
%\date{}

\pagebreak

\begin{center}
\vspace{0.25in} {\large Summary}
\vspace{0.5in}

\parbox{4.8in}{\quad 
%{\em Visual comparisons} of $L_1$ and $L_2$-residuals' plots indicate bad leverage cases
 In several   linear  regression data sets,
%data sets
%of 
$Y (\in R)$ on ${\bf X} (\in R^p),$ {\em visual comparisons} of $L_1$ and $L_2$-residuals' plots indicate bad leverage cases.
The phenomenon is confirmed theoretically by introducing  Location Breakdown Point ({\em LBP})  of a functional $T$:  any point
where  the derivative of $T$'s  Influence Function either takes values at infinities or does not exist.
%$IF'.$ 
Guidelines for the  plots' visual comparisons as diagnostic  are provided.
% for outlier detection.
 The new tools used  include {\em E}-matrix and  suggest  influence  diagnostic {\em RINFIN} which  measures  the  {\em distance} in the {\em derivatives of $L_2$-residuals} at $({\bf x},y)$ from  model $F$ and from gross-error model $F_{\epsilon, {\bf x},y}.$ The larger  {\em RINFIN}$({\bf x},y)$  is, the larger $({\bf x},y)$'s influence in $L_2$-regression residual is. {\em RINFIN}  allows measuring   {\em group influence} of $k$ ${\bf x}$-{\em neighboring} data  cases in a  size $n$ sample
using their average, $(\bar {\bf x}_k,\bar y_k),$  as one case with weight $\epsilon=k/n.$
For   high dimensional, simulated data, 
%and $n$ fixed, 
%indicate that  
%with {\em RINFIN} 
the misclassification proportion of bad leverage cases  in data's   {\em RINFIN}-ordering 
% leverage cases 
decreases to zero as $p$ increases, thus  reconfirming  the {\em blessing of high dimensionality} in the detection of remote clusters.
 The visual  diagnostic and {\em RINFIN} are successful in applications and complement each other.

}
% and for identifying outlying 
%leverage cases.}
% affect the analysis. }
\end{center}

\vspace{.2 in}

\bigskip
{\it Some key words:} \quad Breakdown Point,
% Diagnostic, 
Influence Function, Least Absolute Deviation  Residuals,
Least Squares Residuals,
%Least Absolute Deviation  Residuals,
Leverage, Location Breakdown Point,  Local-Shift-Sensitivity,
Masking, Residual's  Influence Index ({\em RINFIN})
% Wiggling.

\noindent
{\em AMS 2010  subject classifications:}  62-07, 62-09, 62J05, 62F35, 62G35

\pagebreak

\section{\bf Introduction}

%{\em The goal}

\quad  Tukey (1962, p.60) wrote: ``Procedures of diagnosis, and procedures to extract indications rather than extract conclusions, will have to play a large part in the future of data analyses and graphical techniques offer great possibilities in both areas.''
% (Tukey, 1962, p.60). The goal of this work is to provide simple procedures to extract indications for remote cases in least squares($L_2$) linear regression that may affect the  statistical analysis. 
%In the exploratory stage of the analysis, unusual patterns and features of the observations may be revealed by plotting  the residuals against data related quantities. These %quantities are usually the fitted values, the independent variables, the order the  observations were obtained (if it is known) or anything else that is sensible for the particular problem %into consideration (Draper and Smith, 1966). 

In linear
%,  least squares ($L_2$) 
regression of $Y$ on ${\bf X}$  it is often assumed that the data follows probability model $F;  Y \in R, {\bf X} \in R^p.$
However, ``It also happens not infrequently that only part of the data obeys a different model.'' (Hampel {\em et. al.}, 1986). 
Thus, in reality, data may follow  gross-error model $F_{\epsilon,G}=(1-\epsilon)F+\epsilon G$ (Huber, 1964); $G$ is gross-error probability, $0<\epsilon<1.$
The goal of this work is  to provide  {\em simple and fast procedures} for extracting  indications when
% there are 
remote cases  (from
$G$) affect
%ing
 the  statistical analysis in  least squares ($L_2$) regression.  
%Remote cases from $G$ are of particular  interest  because 
%Hampel {\em et. al.} (1986, p. 25, l. -15),  
%``It also happens not infrequently that only part of the data obeys a different model.'' (Hampel {\em et. al.}, 1986). 
These  procedures are particularly  useful for Big Data, when the number of predictors, $p,$ and the sample size, $n,$ increase to infinity;
%  for linear $L_2$-regression  
for $L_2$-regression  $p<n.$

%{\em The motivation}

The initial  motivation
% for this work 
 was provided by the observation,  in  several data sets,  that:  neighboring, remote factor space cases,
$({\bf x},y),$  have  
 least absolute deviation ($L_1$) regression  residuals significantly  larger in size than the corresponding $L_2$-regression   residuals; see, e.g.,  Figures 1 and 2 in section 4. 

 This phenomenon is theoretically confirmed herein using new tools:   {\em E}-matrix and derivatives of the regression coefficients' Influence Functions.
%at $({\bf x},y),$
%that is  not $L_1$-location breakdown point (LBP ),
%The latter allows 
The latter allow 
 calculating  changes in  $L_1$ and $L_2$-regressions residuals  for small perturbations of $({\bf x},y)$ from $F$ and also  from $F_{\epsilon,G}.$
The calculation of $L_1$-residuals' changes  is possible when  $({\bf x},y)$ is {\em not} $L_1$-{\underline L}ocation {\underline B}reakdown {\underline P}oint ({\em LBP}), 
%i.e. the derivative of the Influence Function used does not take value at infinities and 
thus  first order linear approximation of the $L_1$-residual near $({\bf x},y)$  is valid.
{\em LBP} complements the notion of {\em Weight} Breakdown Point (Hampel, 1971). Derivatives of Influence Functions  indicate  a new  influence  diagnostic: {\em RINFIN} (see below).
% In addition, {\em LBP} complements the notion of {\em Weight} Breakdown Point (Hampel, 1971).

A simple graphical method is thus proposed to detect rapidly  in  linear regression data 
%one or more 
remote 
cases, $({\bf x},y),$ affecting drastically $L_2$-regression coefficients.
% when there is one source of ${\bf x}$-gross-error; 
%${\bf x} \in R^p, \ y \in R.$  
Plots of absolute regression residuals against square ${\bf x}$-length provide 
the {\em visual  indications} 
 when 
%Least Absolute Deviation (
$L_1$-residuals' sizes  for  ${\bf x}$-remote  cases are larger, e.g. double in size,  than 
%the sizes 
%those o
%of 
the corresponding  $L_2$-residuals,
thus causing a larger {\em visual gap} in the $L_1$
 plot. A different pattern in the residuals is used to identify other types of outlying cases near a {\em LBP},
% Location Breakdown Point, 
as described in section 4.
%Plots of absolute residuals against square ${\bf x}$-length provide the {\em visual} indications when 
%$L_2$-residuals' sizes  for  ${\bf x}$-remote  cases are significantly smaller than their $L_1$-residuals' sizes,
%causing a larger {\em visual gap} in the $L_1$-plot.
% when there is one source of ${\bf x}$-gross error.

The regression diagnostic, {\em {\underline R}esidual's  {\underline {Inf}}luence {\underline {In}}dex (RINFIN) }for $({\bf x},y),$  is  also introduced
% for each case ${\bf x},y$ 
that measures 
t{\em he distance  in the derivatives of $L_2$-residuals}  when $({\bf x},y)$ follows either probability  $F$ (the model)
or its gross-error mixture 
$F_{\epsilon, {\bf x},y}$, i.e. 
%$F_{\epsilon, {\bf x},y}=
$(1-\epsilon)F+\epsilon \Delta_{{\bf x},y}, \ 0<\epsilon<1, \Delta_{\bf u}$ unit mass at ${\bf u}.$
The larger {\em RINFIN}$({\bf x},y)$ is, the larger $({\bf x},y)$'s influence  in the $L_2$-residual is.
For a group of  remote ${\bf x}$-{\em neighboring} cases from gross-error probability $G,$ with proportion $\epsilon$ in the data,
their group average $(\bar {\bf x}, \bar y)$ is used as one case from 
%$({\bf x},y)$  in 
$F_{\epsilon, \bar{\bf x},\bar y}$ to calculate the
group's influence,  {\em RINFIN}$(\bar {\bf x}, \bar y),$ that depends  also on $\epsilon.$ 
This is an advantage over other methods that use group deletion to determine influence and are exposed {\em i)} to masking from neighboring $G$ cases that remain in the model, {\em ii)}
to a combinatorial explosion due to the very large number of groups to exclude. 
%the influence index.
%$({\bf x},y)$ is replaced in $F_{\epsilon, {\bf x},y}$ by their group averages.

{\em RINFIN} is successful with several known data sets and in simulations, especially   when the dimension $p$ of the data is large.
In simulations with normal mixtures and $n$ fixed, 
%indicate that  
%with {\em RINFIN} 
the misclassification proportion of bad leverage cases  in the   {\em RINFIN} ordering of the data
% leverage cases 
decreases to zero as $p$ increases. The effect of increase in  $p$-values is equivalent to larger  standardized  distance between the means of $F$ and $G$ 
in $F_{\epsilon,G}.$  A similar phenomenon has been observed and confirmed theoretically for mixture densities in a Projection Pursuit cluster detection method (Yatracos, 2013) due to the ``separation'' of the mixtures'components,  measured by their Hellinger's distance, as $p$ increases. 

%To confirm the visual phenomenon, $({\bf x},y)$'s residuals are calculated for  linear regressions with $F$ and $ F_{\epsilon, {\bf x},y}.$ 
%we evaluate  $({\bf x},y)$'s  effect  in the residuals: {\em a)} before entering model $F$ and after, and {\em b)} when $F$'s %gross-error model with $({\bf x},y)$ is ${\bf x}$-perturbed.
 %Thus, residuals of $({\bf x},y)$ are computed 
%with respect to 
%for the linear regression of 
%case
%$({\bf X},Y)$ with model $F$ and  for gross-error models $F_{\epsilon, {\bf x},y}$ and $F_{\epsilon, {\bf x}+{\bf h},y}$ (Huber, %1964); $F_{\epsilon, {\bf x},y}=(1-\epsilon)F+\epsilon \Delta_{{\bf x},y}, \ 0<\epsilon<1, \Delta_{\bf u}$ unit mass at ${\bf u}.$ 
 %and those of the latter with  perturbation model $F_{\epsilon, {\bf x}+{\bf h},y};$ 
%\ {\bf X}, {\bf x}, {\bf h}$ are all in $R^p, \ Y \in R.$
% $||{\bf h}||$ is small, $|| \cdot ||$ is a norm in $R^p.$
 In a  nutshell,
% under mild assumptions,
the justification  for the visual phenomenon and the form of {\em RINFIN} are presented for simple linear regression:\\
% with $|x|$ large and $|h|$ small,
{\em i)} For $\epsilon (>0)$ small in $F_{\epsilon,{\bf x},y}$ residuals are compared, 
\begin{equation}
\label{eq:nutshellLmresbefaftratio}
%0<
\frac{|r_{2,x,y}(x,y)-r_2(x,y)|}{|r_{1,x,y}(x,y)-r_1(x,y)|} \approx C  |r_2(x,y)|,
% C_1  |r_2(x,y)|,
\end{equation}
%{\em ii)} for $\epsilon$ small, $|h|$ small and $|x|$ large,
%\begin{equation}
%\label{eq:nutshellsrrpert}
%\frac {|r_{2, x+h, y}(x+h,y)-r_{2, x,y}(x,y)|}{|r_{1,x+h, y}(x+h,y)-r_{1,x,y}(x,y)|} \approx C_2  |x-EX|;
%\end{equation}
   $r_m$ and $r_{m,x,y}$ are $L_m$-residuals, respectively,  for 
%models 
$F$-regression  and
%that for model 
$F_{\epsilon, x,y}$-regression, $m=1,2; \  C
% \ C_m$
$ is constant,
% \ EX$ is the mean of $X$ under $F,  \ 
$`` \approx''$ denotes approximation.
When $(x,y)$ is gross-error and for $L_1$ and $L_2$ $F$-regressions $r_1(x,y) \approx r_2(x,y),$
with $|r_2(x,y)|>1,$ 
%under model $F_{\epsilon, x,y}$  
from (\ref{eq:nutshellLmresbefaftratio}) it follows for $F_{\epsilon, x,y}$-regression that 
$L_2$ residual of $(x,y)$ is
% more drastically 
reduced more  than its $L_1$ residual,
% reduction 
especially when  $|x|$ is   large (because then  $|r_2(x,y)|$
 is also large).\\
%When  $|x|$ is   large, $|r_2(x,y)|$ is also large and  $L_2$ residual of $(x,y)$ is more drastically reduced  than its $L_1$ residual %in both (\ref{eq:nutshellLmresbefaftratio}) and (\ref{eq:nutshellsrrpert}).
 %When, for example,  $L_1$ and $L_2$ regressions coincide for $F,$  it is clear 
%{\em iii)}
{\em ii)}  $L_2$-{\em residual's  influence index}  of $(x,y)$  from  gross-error  model $F_{\epsilon, x,y}$  
%with respect to  $L_2$-regression
 is
\begin{equation}
\label{eq:ii}
RINFIN(x,y)=
%RINFIN(x,y; \epsilon, L_2)=
\epsilon \cdot \frac {|2r_2( x,y)(x-EX)-\beta_{1,L_2} [(x-EX)^2+Var(X)]|}{Var (X)};
\end{equation}
%$r_2$ is $(x,y)$'s $L_2$-residual with respect to the regression line with
$L_2$-residual ($r_2$), slope ($\beta_{1,L_2}$),  mean ($EX $) and   variance ($Var X$)  are all  under $F.$ 
%all are with respect to model $F$ and $E, \ Var$ denote, respectively, the mean and variance of $X.$ 
%(\ref{eq:ii}) is used to calculate {\em group}  influence of remote neighboring cases because of $\epsilon.$ 
%$\epsilon$ in (\ref{eq:ii}) allows to calculate {\em group}  influence of remote neighboring cases.
%It is based on the 
%{\em RINFIN} is obtained using $x$-derivatives of influence functions, as explained in section 3. The influence from $y$-%derivatives is observed in the residuals' plot; see Remark \ref{r:yresinfl}. 

%Remote ${\bf x}$-cases are studied  since
%Hampel {\em et. al.} (1986, p. 25, l. -15),  
%``It also happens not infrequently that only part of the data obeys a different model.'' (Hampel {\em et. al.}, 1986). 

%{\em Location}  breakdown point ({\em LBP}) 
{\em LBP} of a statistical functional $T$  is motivated and  introduced in section 2 using
${\bf x}$-perturbations of $F_{\epsilon, {\bf x}}.$ 
%gross-error model. 
{\em LBP}  is  a point where the directional or one of the partial derivatives   of $T$'s Influence Function (Hampel, 1971, 1974) 
 either take values at infinities or do not exist. {\em  Local-shift-sensitivity} (Hampel, 1974) cannot replace the 
 derivatives, as explained.
In section 3, regression coefficients' Influence Functions 
%(Hampel, 1971) 
and their derivatives, obtained via $E$-matrices,  are used to show that:
in $F_{\epsilon, {\bf x},y}$-regression, when remote ${\bf x}$-case
becomes slightly more extreme without  reaching  $L_1$  {\em LBP},
%{\em location}-breakdown point,  
 the size of the corresponding 
 $L_2$-residual is drastically reduced whereas the $L_1$-residual is reduced less.
%Regression coefficients' influence functions (Hampel, 1971) 
%and their derivatives, obtained via $E$-matrices  in section 3,  are used to obtain the results.
%{\em  Local-shift-sensitivity} (Hampel, 1974) cannot replace the  derivatives, as explained.
% in section 2.
% Remote ${\bf x}$-cases in models $F_{\epsilon, {\bf x}}$ and $F_{\epsilon, {\bf x},y}$ are studied  since
%Hampel {\em et. al.} (1986, p. 25, l. -15),  
%``It also happens not infrequently that only part of the data obeys a different model.'' (Hampel {\em et. al.}, 1986). 

The graphical method and {\em RINFIN}  are  supported by  applications and simulations in section 4.  Instead of square ${\bf x}$-length on the plot's horizontal axis,
${\bf x}$-length can be used. For some data sets, plotting regression residuals rather than their absolute values may be more informative.
However, for remote gross-error model with small variance, 
%density instead of point-mass   in $R^p,$ with small variance absolute residuals are informative; 
e.g. cases 1-10 in Hawkins-Bradu-Kass (1984) data, absolute residuals are informative.
  
Robust residual plots are accompanied with  confidence ellipsoids.  
Otherwise, visual indications from distorted residuals  lead to inaccuracies.
$L_1$ and $L_2$ residuals and the square length  of dependent variables are used herein because they do not cause
{\em unknown} amount of distortion   in {\em relative} visual distances.
% unlike their robustified versions. 
For example,  in the Stackloss Data  plot
% of standardized Least Median of Squares (LMS) residual against a robust distance 
% for Stack Loss Data in
(Rousseuw and van Zomeren,1990, p. 636, Figure 3)
%: {\em a)}
% robustified ${\bf x}$-distance of case 4  is comparable to the smallest robustified ${\bf x}$-distances whereas  square length of %case 4  in Figure 1 herein  is larger than half (or more) square lengths in the data.
% {\em b)} the {\em absolute} standardized LMS  residuals of cases 1, 3,  4 and 21
%are  {\em all} comparable, whereas in  Figure 1, the absolute residual of case 21 is almost double than those of cases 1,  3 and %larger than that of case 4; in  our Figure 2,  the corresponding  sizes are smaller but relative sizes remain comparable.
relative sizes  of the {\em absolute} standardized Least Median of Squares (LMS)  residuals of cases 1, 3,  4 and 21 differ  from those  in the $L_2$-absolute residuals in Figure 2 herein.

In multiple regression, with 
%the response and the independent variables 
observations from  
%probability model 
$F,$ a case $({\bf x}, y) $ with  
%remote
 factor space component, ${\bf x},$
%-component  
far away from the bulk of $F$'s  factor space
% in $F$ 
  is
called {\em leverage case}  (Rousseeuw and Leroy, 1987,
% p. 6, l. -1 to -3; 
Huber, 1997
%, p. 60, l. 1
). A ``good'' leverage case 
%$({\bf x}, y) $
 is either  near or on the regression hyperplane
determined by $F.$
% determined by the majority of the data.
 A ``bad'' leverage case  forces the $F$-hyperplane to change drastically when ${\bf x}$ becomes more remote.
% with respect to  because of  its  outlying $y$-response. 
%Bad leverage cases are not easily identified.
The suggested  comparisons of $L_1$ and $L_2$ residuals' plots  and 
data's {\em RINFIN} values  reveal ``bad'' leverage cases.

 The abundance of high dimensional data sets from various fields, with $p$ and $n$  both large,
% and in particular $n<p,$ 
created the need for new methods to detect outliers/influential cases affecting linear regression analysis. High dimensional influence measures have been proposed among others by Alphons {\em et al.} (2013), Zhao {\em et al.} (2013, 2016) and  references there in. She and Owen (2011) identify influential cases using nonconvex  penalized likelihood.
 Influence Function in outlier detection has been used by Campbell (1978) and Boente {\em et al.} (2002).
 The influence of observations in estimates' values  has been also studied by several authors, among others by  Cook (1977), Cook and Weisberg (1980), Ruppert and Carroll (1980), Carroll and Ruppert (1985),  Hampel (1985),  Hampel {\em et. al.} (1986),  Ronchetti (1987), Rousseeuw and van Zomeren (1990), Ellis and Morgenthaler (1992),  Bradu (1997),  Flores (2015) and Genton and Hall (2016).
%Campbell (1978) and Boente {\em et al.} (2002) used Influence Function in outlier detection.

In Genton and Ruiz-Gazen (2010)
%characterized 
an {\em observation} is
influential ``whenever a change in its value leads to a radical change in the estimate'' 
and  the {\em hair-plot} is introduced to identify it.
 %(p. 808, l. -8 to -7). 
 Two influence  measures  are proposed  using  partial derivative  of the {\em estimate}: {\em a)} the local, 
 with a small perturbation in one coordinate of the observation, and {\em b)} the global,
% (or asymptotic influence), 
using the most extreme contamination for each  coordinate. 
%Both measures are computed for quadratic forms  (of the observations) and for their ratios.
%The obtained local influence measures take finite values and the global measure may become infinite.
 Differences in 
%the approach in  
our work include:   {\em i)}
% our main goal  is  to
%we  identify visually, using $L_1$-regression residuals,
% linear regression 
leverage cases affecting drastically $L_2$-regression residuals
are visually identified {\em combining information}  from $L_1$ and $L_2$ 
residuals' plots,
{\em ii)}  the  
%directional and partial 
derivative of the estimate's 
%a functional's 
 influence function  is used instead of
% that of 
the estimate's derivative,
% to approximate locally  functionals' values,
{\em iii)} {\em RINFIN}  measures distance in 
residuals' derivatives and can be used to evaluate group  influence of neighboring cases.
% of cases in the residuals' change, in particular group influence, 
%  {\em iv)} one additional goal is to show the usefulness in  combining visual information from functionals and estimates  with high %and low breakdown points.
 
Work has   been done  to identify  ``bad'' leverage cases  using $L_1$ residuals.
Barrodale (1968) compared $L_1$ and $L_2$
residuals for regression function
$\sum_{j=1}^p a_j \phi_j({\bf x})$
using tables for different ${\bf x}$-values;  ${\bf x}  \epsilon R^d, \phi_j$ is known,
$a_j$ is an  unknown coefficient, $j=1,...,p.$
%and examined the effect of ``wild'' cases;
%$ x \epsilon R^d, \phi_j$ is known, $a_j$ is an  unknown coefficient, $j=1,...,p.$
%Bloomfield and Steiger (1983) present a unified treatment of the role of the least absolute deviation techniques
%in linear regression and other domains including theory, applications and algorithms.
Barnett and Lewis (1984) present the absolute residuals
as a tool in outlier detection.
% along with other related topics and several references.
%In a Monte Carlo study of the estimation of a location parameter Hampel (1985) used the notion of (weight) breakdown point
%to explain the behavior of several outlier tests and the masking effect.
Narula and Wellington (1985) looked for observations that do and do not
affect the analysis in $L_1$-regression using the residuals.
%Dielman (1986) observed in Monte Carlo simulations that $L_1$ regression
%forecasts seem more efficient than $L_2$ forecasts in heavy tail distributions.
%Ronchetti (1987) presents a review of the robustness approach based on influence functions with applications to linear models. %The influence function is used to assess, among others, techniques based on the $L_1$ distance and to motivate new robust %procedures.
Ellis and Morgenthaler (1992) and Bradu (1997) examined the performance of the $L_1$ regression estimator facing outliers
in the {\em response} variable.
%More research 
Additional results on $L_1$-regression and outliers may be found in Dodge(1987).

%More work has been  devoted recently on 
%the study of 
Recent  results combine also   information from $L_1$ and $L_2$ regression.
 %have been obtained on the relation between the $L_1$ and the $L_2$-residuals and the so-derived  information
%that can be obtained. 
%In a survey on linear regression methods, 
Giloni and Padberg (2002) presented a lower bound on total sum of absolute $L_1$-residuals using the total sum of  squared $L_2$-residuals. 
Flores (2015) studied for a {\em particular} regression model  the behavior  of 
%the
 $L_1$-estimates   
%by comparison  
by comparing them with  $L_2$ -estimates, 
%and concludes for this particular model that ``the $L_1$-estimator mitigates the effect of the largest outliers.'' This is confirmed  herein for 
 and introduced  {\em leverage constants} for a design matrix
to determine whether leverage cases are good or bad.

Proofs are in the Appendix.

\section{\bf Location Breakdown Point ({\em LBP})}

%Some preliminary notions are presented that may be found in
%Barnett and Lewis(1984, p. 68-74),
%Hampel, Ronchetti, Rousseeuw and Stahel(1986), Huber(1981, 1997)  and in
%Rousseeuw and Leroy(1987, p. 186-188).

\quad Hampel (1971) introduced the influence function, $IF({\bf x}; T, F),$ 
of a functional  $T$ at 
%a cumulative distribution (or
 probability  $F,$ 
% evaluated at $x (\inR^p$) 
%as 
% (Hampel, 1971)
% is
%introduced in Hampel (1971),
% and is defined as
\begin{equation}
\label{eq:inffunc}
IF({\bf x}; T, F) =
\lim_{\epsilon \rightarrow 0}\frac{T[(1-\epsilon)F+\epsilon \Delta_{{\bf x}}]-T(F)}
{\epsilon}, 
\end{equation}
%at  points 
%$x(\epsilonR^p)$ for which 
when this limit exists; ${\bf x} (\in R^p), \ \Delta_{{\bf x}}$ is the probability
distribution that puts all its mass at the point ${\bf x}, \ 0<\epsilon<1.$  
%The influence function

$IF({\bf x}; T, F)$ 
determines
% measure of  
the ``bias'' in the value of $T$ at $F$   due to an $\epsilon$-perturbation  of $F$ with $\Delta_{{\bf x}}:$
%adding a proportion $\epsilon$ of
%outliers at
% a point 
%${\bf x}$ via  the linear approximation:
%of the estimator:
\begin{equation}
\label{eq:firstorder}
T[(1-\epsilon)F+\epsilon \Delta_{{\bf x}}] - T(F) \approx  \epsilon IF({\bf x}; T,F).
\end{equation}

%For the sample version of influence $T(\hat F_n)$ is used rather than $T(F);$
%$\hat F_n$ is the empirical distribution of the sample $x_1,...,x_n.$
%Using the previous notation and a fraction $\epsilon=\frac{m}{n}$
%of outliers the approximation
%$$T(x_1,...,x_{n-m}, x,...,x) \approx T(F_n)+ \frac{m}{n}IF(x; T,F)$$
%is obtained.
\bef (Hampel, 1971) The weight  breakdown point 
%$\epsilon_n$ 
is the upper bound
%$\frac{m}{n}$ 
on $\epsilon$ for which
 linear approximation (\ref{eq:firstorder}) can be used.
\enf
\quad Discussing further concepts related to the influence function,  Hampel (1974, p. 389) introduced
% the value of
 {\em  local-shift-sensitivity},
\begin{equation}
\label{eq:LSS}
\lambda^*=\mbox{ sup}_{{\bf x}
\neq {\bf y}}
\frac { |IF({\bf x};T,F)-IF({\bf y};T,F)|}{||{\bf x}-{\bf y} ||}, 
\end{equation}
%and the derivative of the influence function
as ``a measure for the {\em worst} (approximate) effect of wiggling the
observations''; 
%${\bf x} \in R, y \in R, \
$ || \cdot || $ is a Euclidean distance in $R^p.$
%adding that 
%without connecting it to the derivative of the Influence Function.
%It was mentioned that
%``one has to keep in mind that it refers only to standardized local changes of the value of the estimator, so that even an infinite %value of $\lambda^*$  may refer only to a very limited actual change''.  

Unlike the extensive use of the {\em weight} breakdown point, 
%its empirical counterpart  and its extensions,
{\em  local-shift-sensitivity} was never fully 
exploited. One reason is that,  in reality,  it is  a ``global'' measure as supremum over all  ${\bf x}, {\bf y}.$ Thus, 
$\lambda^*$  cannot be used to study $T$'s 
%change of value 
bias for ${\bf x}$'s small perturbation in the $\epsilon$-mixture, from ${\bf x}$ to ${\bf x}+{\bf h}, \ ||{\bf h}||$  small, 
\begin{equation}
\label{eq:biassmallchanges}
T[(1-\epsilon)F+\epsilon \Delta_{{\bf x}+{\bf h}}]-T[(1-\epsilon)F+\epsilon \Delta_{\bf x}].
%= \epsilon h \frac{T[(1-\epsilon)F+\epsilon \Delta_{x+h}]-T[(1-\epsilon)F+\epsilon \Delta_x]}{\epsilon h},
\end{equation}  

%The local-shift-sensitivity is a {\em global}  measure of ``wiggling'' and  cannot be used to evaluate the difference
%$$T[(1-\epsilon)F+\epsilon \Delta_{x+\eta}]-T[(1-\epsilon)F+\epsilon \Delta_x]$$
%for different remote $x$-values and for small $\eta.$ This can be accomplished using instead the derivative of the influence %function. However, it is seen herein that  a ``limited change''   is   informative  in plots' comparisons for  $L_1$ and $L_2$ residuals.

%However, this limited actual change is informative
%in $L_1$-residual plots.

Rousseeuw and Leroy (1987)
presented a physical analogy to  the notion of {\em weight}  breakdown point.
A beam is fixed at one end and, at point ${\bf x}$ on the beam, a stone with weight $\epsilon$ is attached.
%to which a stone with weight $\epsilon$ is attached,  exercising its force at the point ${\bf x}.$ 
For small weights, the ``deformation'' (i.e., the bias) (\ref{eq:firstorder})
%(denoted w.l.o.g. by)  $T(x)$ 
%$$T[(1-\epsilon)F+\epsilon \Delta_{\bf x}] -T(F),$$ 
of the beam is linear
in $\epsilon$ and one can predict the weight's effect.
% linear approximation  (2).
As soon as $\epsilon$ takes value larger than
the ``breakdown value'' (that depends on the location ${\bf x}$), (\ref{eq:firstorder})  cannot be used.

For the physical analogue of  {\em location}  breakdown,  a  sufficiently long beam is used  and 
%we let  the (same) 
weight
$\epsilon$ ``travels'' at different   ${\bf x}$-locations
% along a line, 
 far away
from the fixed end of the beam. 
%For the deformation 
%$T(x)$
There is a location ${\bf x}_{0,\epsilon}$
that  makes the beam ``break''.
The beam will break also with a small perturbation from ${\bf x}_{0,\epsilon}-{\bf h}$ to ${\bf x}_{0,\epsilon}, \ ||{\bf h}||$ small.
This is the reason we study ${\bf h}$-perturbations (\ref{eq:biassmallchanges})  for remote ${\bf x}$'s.
%and  linear approximation (\ref{eq:firstorder}) 
%to $T[(1-\epsilon)F+\epsilon \Delta_{x_{0,\epsilon}}] -T[(1-\epsilon)F+\epsilon \Delta_x] $
%will not be valid
%; $x$ is at a shorter distance from the fixed end
%of the beam than 
%at ${\bf x}={\bf x}_{0,\epsilon}.$

When $F$ is defined on the real line, 
to express the physical analogue of {\em location} breakdown with the derivative of the influence function we 
 evaluate (\ref{eq:biassmallchanges}) at    neighboring points $x, \  x+h, \  x \in R, \  h  \in R, |h|$  small.
%The key Lemma in this work follows.
\bel 
%(A new tool)
\label{l:key}
 \begin{equation}
\label{eq:limit1}
\lim_{h \rightarrow 0}  \lim_{
%h \rightarrow 0 
\epsilon \rightarrow 0}
\frac{T[(1-\epsilon)F+\epsilon \Delta_{x+h}]
-T[(1-\epsilon)F+\epsilon \Delta_x]}
{\epsilon h }=IF'(x;T,F).
%=\infty.
\end{equation}
\enl

%The derivative of the influence function
$IF'(x;T,F)$  is  used to approximate (\ref{eq:biassmallchanges}) for small $\epsilon, \  |h|:$   
%difference (\ref{eq:biassmallchanges}),
\begin{equation}
\label{eq:addedvalue}
T[(1-\epsilon)F+\epsilon \Delta_{x+h}]-T[(1-\epsilon)F+\epsilon \Delta_x] \approx \epsilon h IF'(x;T,F);
\end{equation}
(\ref{eq:addedvalue}) is {\em the  tool} used to 
%obtain  approximations for
approximate   $L_1$ and  $L_2$ residuals and determine 
group influence.

In simple  linear $L_1$-regression, derivatives of $IF(x,y; T, F)$ are constant  
%with case $(x,y)$ instead of $x$  but $y$ fixed,
% at cases where
where the residual does not vanish;
%$IF'(x,y; T, F)$ is constant; 
$T$ is any of the  regression coefficients.
% $IF'(x; T, F)$ is constant.
% and $IF'(x; T, F)$ is constant for most $x$-values. 
As 
%the case 
$(x,y)$ becomes more remote in $x,$
eventually there is a   change of the $L_1$-regression coefficients at $(x+h,y),$ 
%and continuity is lost,
the $L_1$- residual vanishes, 
 %However, in location points where
%$|IF'(x+h; T, F)|$ 
 derivatives of $IF(x,y; T, F)$ takes values infinities and   (\ref{eq:addedvalue}) is {\em not} valid.
%there is a substantial change of the $L_1$-regression coefficients.

%For  $L_1$-regression, possible consequences of adding or dropping a single data point or some subset of data are ``either
%no effect, or passing from a unique fit to another one, or reduce or increase the extent of non-uniqueness of fit that is already
%present''(Bloomfield and Steiger, 1983, p. 57). Thus, if $(x+h,y)$ replaces $(x,y)$ in $F_{\epsilon,x},$   eventually a  different fit %occurs for some $h$  with  loss of continuity in the regression coefficients and the  $L_1$-residuals' values.

This observation  motivates the definition of {\em location} breakdown point ({\em LBP})  where the derivative of the influence function takes infinite values. In $L_1$ and $ L_2$ linear regressions partial derivatives of the coefficients influence functions  exist and, in addition, one  remote  coordinate in the factor space is enough to reach {\em LBP}.
%a {\em location} breakdown point.
Thus, in the definition of {\em LBP} for $T$ partial  derivatives are used  instead of a directional derivative.

%where (\ref{eq:addedvalue}) is not valid.
%regression coefficients may not change substantially  but the residual will have larger absolute 
%value when $x$ is a leverage case. In linear $L_2$-regression  $IF'(x; T, F)$ is linear in $x$ and both the regression coefficients; %and the residuals are affected when $x$ is a leverage case.
%It point $x \in R$ for which $|IF'(x; T, F)|$ is infinity 
%(\ref{eq:locbreakdownG})
% is a location breakdown point. This motivates the following definition for $x \inR^p, \ k\ge 1.$

\bef
\label{f:locbreakdownG}
Let $T$ be a functional defined on probabilities in $R^p,$ with real values,  $p \ge 1.$
 Then,
% a point 
${\bf x} \in R^p$
is Location Breakdown Point (LBP) if 
%{a)}
 there is $j \le p:$
%any partial derivative of the influence function $IF(y; T,F)$ takes value at infinities when $y=x.$
\begin{equation}
\label{eq:locbreakdownG}
% |\frac{\partial IF({\bf x};T,F)}{\partial x_j}|=\infty; 
 |\frac{\partial }{\partial x_j} IF({\bf x};T,F)|=\infty  \mbox{ or does not exist;}
\end{equation}
$x_j$ is ${\bf x}$'s $j$-th coordinate, $F$ is probability.
%{b)} ${\bf x}$ is far away from the bulk of the data.
\enf
%\ber
%Since $F$ is usually unknown, location breakdown points may be estimated using $\frac{\partial }{\partial x_j} IF({\bf x};T, \hat %F_n);
%\hat F_n$ is the empirical distribution, $j \le k.$
%\enr
% However, this limited actual change is informative in $L_1$-residual plots.
%The  x-values for which $IF'(x; T, F)$ does not exist  constitute the location breakdown points. 
%Except of
%the physical analogy this notion is motivated by the behavior of
%estimates in $L_1$ regression; see next section.
%Examples follow with  $F$  defined on measurable subsets of $R$ and $R^2.$
%calculating $IF'(x;T,F),  \ x \in R,$ evaluated at the (population) distribution $F.$ 
%the sample version is obtained replacing F by the empirical cumulative distribution function $\hat F_n.$ 
%The von Mises differentiation technique
%as described in Reeds (1976) is followed to obtain the influence functions.
% Since $F$ is usually unknown, location breakdown points may be estimated using $IF'(x;T, \hat F_n).$
\beex
{\em  Let   $F$ be a probability on the real line,
%\ F'(x)=f(x), EX^2<\infty,$
% $T_1(\hat F_n)=\mbox{ median}(X_1,...,X_n)$
$T_1(F)$ is the median of $F,$ 
%and $T_2(\hat F_n)=\bar X.$
$T_2(F)$ is the mean of $F$ and their influence functions are:
%It can be seen that
$$IF(x;T_1,F)=\frac{sign[x-T_1(F)]}{2f[T_1(F)]}, \hspace{5ex}
IF(x;T_2,F)=x-T_2(F).$$
From  (\ref{eq:locbreakdownG}), there are no 
%location breakdown points 
{\em LBPs} on the real line  for the mean, $T_2,$
% of the observations
 but for the median, $T_1,$ 
%of the observations
its value
is the only {\em LBP}.}
%{\em location} breakdown point.}
\enex

\beex
{\em Consider a simple linear regression model, $Y=\beta_0+ \beta_1 X + e,$
with
% the distribution of
 error $e$ 
%symmetric around zero and having
having mean zero and finite second
moment, $F$ is the joint distribution
of $(X,Y) \mbox{ and } f_{Y|X}$ is the conditional density of $Y$ given $X,$
% \ X \in R, Y \in R.$
\begin{equation}
\label{eq:ftilde}
\tilde f_{Y|X}(x)= f_{Y|X}[\beta_{0,L_1}(F)+ \beta_{1,L_1}(F)x|x].
\end{equation}

The influence functions for
the $L_2$ 
%estimates 
-parameters $\beta_{0,L_2}(F), \  \beta_{1,L_2}(F),$ obtained at  $F$ 
are 
%$\hat \beta_{0,L_2}=\beta_{0,L_2}(\hat F_n), \hat \beta_{1,L_2}=
%\beta_{1, L_2}(\hat F_n)$ are given in terms of $F$ by:
\begin{equation}
\label{eq:beta0L2}
IF(x,y;\beta_{0,L_2}(F), F)=[y-\beta_{0,L_2}(F)-\beta_{1,L_2}(F)x] \frac {EX^2-xEX} {Var(X)}, 
\end{equation} 
\begin{equation}
\label{eq:beta1L2}
IF(x,y; \beta_{1,L_2}(F), F)=[y-\beta_{0,L_2}(F)-\beta_{1,L_2}(F)x] \frac {x-EX}{Var(X)};
\end{equation} 
$EU$ and $Var(U)$ denote, respectively, $U$'s mean and variance.
The derivatives of influence functions (\ref{eq:beta0L2}), (\ref{eq:beta1L2})  do not satisfy (\ref{eq:locbreakdownG}) for $x \in R, \ y \in R, $ thus there are no {\em LBPs}.
% location breakdown points. 
%Note that the $y$-derivatives of the influence functions  are respectively $(EX^2-xEX)/Var(X)$ and $(x-EX)/Var(X)$ thus do not %contain information about $r(x,y).$

The influence functions for the $L_1$-parameters $\beta_{0,L_1}(F), \  \beta_{1,L_1}(F),$ obtained at  $F$
are
%$\hat \beta_{0,L_1}=\beta_{0,L_1}(\hat F_n), \hat \beta_{1,L_1}=
%\beta_{1,L_1}(\hat F_n)$ are given in terms of $F$ using the notation
%\begin{equation}
%\label{eq:ftilde}
%\tilde f_{Y|X}= f_{Y|X}[\beta_{0,L_1}(F)+ \beta_{1,L_1}(F)x|x]
%\end{equation}
% by:
\begin{equation}
\label{eq:beta0L1}
IF(x,y; \beta_{0,L_1}(F), F)=\frac{
sign[y-\beta_{0,L_1}(F)-\beta_{1,L_1}(F)x]}{2}
\frac {EX^2 \tilde f_{Y|X}(X)
-xEX \tilde f_{Y|X}(X)}
{E \tilde f_{Y|X}(X) EX^2 \tilde f_{Y|X}(X)
 -(EX \tilde f_{Y|X}(X))^2},
\end{equation}
\begin{equation}
\label{eq:beta1L1}
IF(x,y;\beta_{1,L_1}(F), F)=\frac{sign[y-\beta_{0,L_1}(F)-\beta_{1,L_1}(F)x]}{2}
\frac {xE \tilde f_{Y|X}(X)
-EX \tilde f_{Y|X}(X)}
{E \tilde f_{Y|X}(X) EX^2 \tilde f_{Y|X}(X)
 -(EX \tilde f_{Y|X}(X))^2};
\end{equation}
%the constant $\tilde f_{Y|X}$ is 
%\begin{equation}
%\label{eq:ftilde}
%\tilde f_{Y|X}= f_{Y|X}[\beta_{0,L_1}(F)+ \beta_{1,L_1}(F)x|x].
%\end{equation}

From (\ref{eq:locbreakdownG}), 
{\em LBPs}
% {\em location}  breakdown points 
in $L_1$-regression
 are all  $x,y$ satisfying the relation
$y=\beta_{0,L_1}(F)+\beta_{1,L_1}(F)x.$ 
%The $y$-derivatives of the influence functions either vanish or take values at infinities.
}
\enex

\ber
\label{r:noinfoiny}
The $y$-derivatives of $L_2$-influence functions (\ref{eq:beta0L2}), (\ref{eq:beta1L2}) are, respectively,
$(EX^2-xEX)/Var(X)$ and $(x-EX)/Var(X);$ those of $L_1$-influence functions (\ref{eq:beta0L1}), (\ref{eq:beta1L1}) 
either vanish or take values at infinities. 
%Thus, $y$-derivatives of influence functions do not provide information for $r(x,y)$
%and often we use $y$ as fixed.
\enr

%\section{\bf $IF'({\bf x};T,F),$ Residuals,  Leverage Cases,  Regression Residuals}
% Leverage Cases  and Location Breakdown Points}
\section{Influence, Residuals,  Leverage Cases, {\em RINFIN}}

 {\em MULTIPLE REGRESSION MODEL}

Let $({\bf X},Y)$ follow probability model $F$ in $R^{p+1},$
\begin{equation}
\label{eq:mulregmod}
Y= {\bf X}^T \beta+ e;
\end{equation}
 ${\bf  X}=(1, X_1, ...X_p)^T$
is the 
%vector of 
independent variable,  $Y$ is the response, $\beta=(\beta_0,..,\beta_p)^T.$ 
%$\sigma_i^2$ is the variance of $X_i.$
%is the $(p+1)$x$1$ unknown vector parameter of interest.  

{\em The Model Assumptions:}\\
 (${\cal A}1$)
The error, $e,$ is symmetric around zero and has finite second moment.\\
%{\em b)} {\em w.l.o.g.}  Var X=1; to avoid one extra denominator in the influence functions.\\
%symmetric error $e$  with finite second moment, let $F$ be the cumulative distribution of $(X,Y).$ 
(${\cal A}2$) $X_1, \ldots, X_p$ are independent random variables. \\
(${\cal A}3$) Case $({\bf x},y)$ is mixed with cases from model $F$  with probability $\epsilon$ (model $F_{\epsilon, {\bf x},y}$).

 Let $({\bf x}+{\bf h},y), ({\bf x},y+h) $ be  small perturbations of
% more remote  from the bulk of the data than 
$({\bf x},y).$ 
%and assume it replaces it. 
The goal is to compare the $({\bf x},y)$- residual  changes  in  $L_1$ and in  $L_2$ regressions:\\
{\em i)}  before $({\bf x},y)$ enters model $F$ and after, i.e.,  under $F_{\epsilon, {\bf x}, y},$   \\ {\em ii)} when $({\bf x}+{\bf h},y)$ replaces $({\bf x},y)$ in the $\epsilon$-mixture,  i.e.,  under $F_{\epsilon, {\bf x}, y}$ and  $F_{\epsilon, {\bf x}+{\bf h}, y}$ and\\ 
{\em iii)} when $({\bf x},y+h)$ replaces $({\bf x},y)$  in the $\epsilon$-mixture,  i.e.,  under $F_{\epsilon, {\bf x}, y}$ and  $F_{\epsilon, {\bf x}, y+h}.$
% To examine the effect of leverage we assume $y$ to be fixed and study  ${\bf x}$-perturbations.
%Comparison of residuals under $F_{\epsilon, {\bf x},y}$ and $F_{\epsilon, {\bf x},y+h}$ indicate large influence of
%cases at the data's extremes in plot of residuals; see Remark \ref{r:yresinfl}.  

Let  ${\bf x}$  become  more extreme in the $i$-th coordinate, $x_i+h, \ |h|$ small; denote by ${\bf x}_{i,h}$ this perturbation of ${\bf x},$
\begin{equation}
\label{eq:xperturb}
{\bf x}_{i,h}={\bf x}+(0,\ldots, h,\ldots,0).
\end{equation} 

 The $j$-th regression coefficients obtained by $L_m$-minimization, respectively,   at
 models $F_{\epsilon, {\bf x},y}$ 
%$(1-\epsilon)F+\epsilon \Delta_{{\bf x},y}$ 
and $F$
are:  
%denoted by
%$$\beta_{j, L_m,{\bf x}}(\epsilon)= \beta_{j, L_m}([(1-\epsilon)F+ \epsilon \Delta_{{\bf x},y}, ]), \ \beta_{j, L_m}=\beta_{j, %L_m}([F]),  \  j=0,1, ..., p,$$
\begin{equation}
\label{eq:coef}
\beta_{j, L_m,{\bf x}}= \beta_{j, L_m}([F_{\epsilon, {\bf x},y}]), \ \beta_{j, L_m}=\beta_{j, L_m}([F]),  \  j=0,1, ..., p,
\end{equation}
% \ m=1,2,$$
\begin{equation}
\label{eq:coefvec}
\beta_{L_m,{\bf x}}=(\beta_{0, L_m,{\bf x}},...,\beta_{p, L_m,{\bf x}})^T, \hspace{3ex} \beta_{L_m}=(\beta_{0, L_m},...,\beta_{p, L_m})^T;
\end{equation}
% Denote the $L_m$- residual by $r_m({\bf u},v;\epsilon)=v-{\bf u}^T\beta_{L_m,{\bf u}}(\epsilon)$
%and  expected value by $E; \  m=1,2.$
denote 
%for fixed $v$ 
the $L_m$- residuals  for models $F_{\epsilon,{\bf u},v}$  and $F,$ respectively, 
\begin{equation}
\label{eq:resnot}
r_{m, {\bf u}}=r_m({\bf u},v;F_{\epsilon,{\bf u},v})=v-\beta_{L_m,{\bf u}}^T{\bf u}, \hspace{3ex} 
r_m=r_{m}({\bf u},v)=v-\beta_{L_m}^T{\bf u}, \hspace{3ex} m=1,2.
\end{equation} 

%{\em Notation:} 
When indices of $\beta$'s and $r$  include at least one among ${\bf x},{\bf x}_{i,h}, {\bf u}, y+h,$ they are determined from a gross-error model.  Only ${\bf x}$ is used  at $\beta_{j, L_m,{\bf x}}$ and only  ${\bf u}$ is used at $r_{m, {\bf u}}$ because of interest
in factor space perturbations and to avoid increasing  the number of indices.  The influence function of $\beta_{j,L_m}$ is evaluated at $({\bf x},y)$ for $F,$  thus use 
%the $j$-th regression coefficient and its partial derivative in $x_i,$
%$IF_{j,L_m}$ of $\beta_{j,L_m}(0)=\beta_{j,L_2}([F]), \  j=0,1, \ldots, p, $ by
\begin{equation}
\label{eq:IFnotation}
IF_{j,L_m}=IF({\bf x},y;\beta_{j,L_m},F), 
 \hspace{5ex} IF'_{v,j,L_m}=\frac {\partial IF({\bf x},y;\beta_{j,L_m},F)}{\partial v}, \hspace{3ex} v=y, x_i,
\end{equation}
i.e., in words, $ IF'_{v,j,L_m}$ is the derivative of  $IF_{j,L_m}$ with respect to $v,$
$i=1, \ldots, p, \  j=0,1, \ldots, p, \  m=1,2.$

 Influence functions of $L_m$ regression coefficients are
solutions of  the  equations:
\begin{equation}
\label{eq:syst1}
 IF_{0,L_m}  + IF_{1,L_m} EX_1+...+ IF_{p,L_m}EX_{p}= \tilde  r_m({\bf x},y),
\end{equation}
%$$  IF_{0,L_m} EX_i + ...+ IF_{p,L_m}(0) EX_iX_j+...+ IF_{p,L_m}  EX_iX_{p}=x_i \tilde r_m({\bf x},y)$$
\begin{equation}
\label{eq:syst2}
 IF_{0,L_m} EX_i + ...+ IF_{p,L_m} EX_iX_j+...+ IF_{p,L_m}  EX_iX_{p}=x_i \tilde r_m({\bf x},y),  \mbox{ } i=1, \ldots ,p, \  m=1,2,
%=\int \int r_2({\bf u},v;0)x_i  d(\Delta_{{\bf x},y}-F)({\bf u},v),  \mbox{ } i=1, \ldots ,p.
\end{equation}
\begin{equation}
\label{eq:rtilde}
\mbox{with } \hspace{6ex} \tilde r_1({\bf x},y)=\frac{sign[r_1({\bf x},y)]}{2 \tilde f_{Y|{\bf X}}}, \hspace{6ex} \tilde r_2({\bf x},y)=r_2({\bf x},y);
\end{equation}
from the symmetry of $e$ in  assumption (${\cal A}1$),   $\tilde f_{Y|{\bf X}}$ 
is  the common value
\begin{equation}
\label{eq:ftildmult}
\tilde f_{Y|{\bf X}}=f_{Y|{\bf X}}({\bf x})=f_{Y|{\bf X}}[\beta_{0,L_1}+ \beta_{1, L_1} x_1+\ldots+\beta_{p, L_1} x_p|{\bf x}].
\end{equation} 
 
\pagebreak

{\em E-MATRIX AND ITS COFACTORS}

%The influence functions for the regression coefficients are solutions of the system  (\ref{eq:syst1}), (\ref{eq:syst2}) .
 Under assumption $({\cal A}2),$ 
the coefficients   in the system of equations  (\ref{eq:syst1}), (\ref{eq:syst2}) 
form a special type of matrix we call $E_p$-matrix; $p$ is the covariates' dimension. As an illustration,
for real numbers $a,b,c,A,B,C,$ 
$$ E_4=
 \left  ( 
\begin{array}{cccc}
1 & a & b & c\\
a & A & ab & ac\\
b & ba & B & bc\\
c & ca & cb & C 
\end {array} 
\right ).
$$
%\end{equation}
 For $E_4,$  the corresponding  linear regression model with independent covariates   $X_1, \ X_2, \ X_3$ 
 provides  
%entries   
$a=EX_1, \ b=EX_2, \ c=EX_3$ and $A=EX_1^2, \ B=EX_2^2,  \ C=EX_3^2.$
%These entries and the pattern  appear when calculating influence functions for regression coefficients' functionals
%in (\ref{eq:systl2}), (\ref{eq:systl1})  and the covariates are independent.
%It follows that $V$ and $c_{i,i}, i=1,\ldots,4$ are all positive. 

\bef $E_n$-matrix  with real entries has form:
\begin{equation}
\label{m:Eleni}
 E_n=
 \left  ( 
\begin{array}{cccc}
1 & a_1 & a_2 \ldots  & a_n\\
a_1 & A_1 & a_1 a_2  \ldots & a_1 a_n\\
a_2 & a_2a_1 & A_2 \ldots & a_2a_n\\
\ldots \\
a_n & a_na_1 & a_na_2 \ldots  & A_n 
\end {array} 
\right ).
\end{equation} 
\enf

{\em Notation:} $E_{n,-k}$ denotes the matrix obtained from $E_n$ by deleting its $k$-th column and 
$k$-th row, $2 \le k \le n+1.$

{\em Property of $E_n$-matrix:} Deleting the $k$-th row and the $k$-th column of $E_n$-matrix, the obtained matrix
$E_{n,-k}$ is  $E_{n-1}$ matrix formed by $\{1,a_1,\ldots,a_n\}-\{a_{k-1}\}, \ 2 \le k \le n+1.$

The cofactors of $E_n$-matrix are needed to solve (\ref{eq:syst1}), (\ref{eq:syst2}).

\bep 
\label{p:Elenimatrix}
a) The determinant of   $E_n$-matrix (\ref{m:Eleni}) 
is
\begin{equation}
\label{eq:ElenimatrixDet}
|E_n|=\Pi_{m=1}^n (A_m-a_m^2).
\end{equation}
b) Let $C_{i+1,j+1}$ be the cofactor of element $(i+1,j+1)$ in $E_n.$
Then, its determinant 
\begin{equation}
\label{eq:cof1}
C_{i+1,j+1}=0, \mbox{ if } i>0, j>0, i \neq j, \hspace{8ex}  C_{1,j+1}=-a_j \Pi_{k \neq j} (A_k-a_k^2).
\end{equation}
 \begin{equation}
\label{eq:cof2}
C_{i+1,1}=-a_i \Pi_{j \neq i} (A_j-a_j^2)       , \mbox{ if } i>0, \hspace{8ex}  C_{1,1}=|E_n|+\sum_{k=1}^n a_k^2 |E_{n, -k}|.
\end{equation}

\enp

%\pagebreak

{\em $L_m$-REGRESSION INFLUENCE FUNCTIONS, m=1, 2}

\bep
\label{p:IF}
For regression model (\ref{eq:mulregmod}) with assumptions (${\cal A}1$)-(${\cal A}3$), $r_1({\bf x},y) \neq 0,$ and  $\tilde r_m$  in (\ref{eq:rtilde}),
%and $\sigma_i^2$  the variance of $X_i$ 
the influence functions
of $L_m$-regression coefficients, $m=1,2,$ are:
\begin{equation}
\label{eq:IFALL}
IF_{0,L_m}=\tilde r_m [1-p+\sum_{j=1}^p \frac{EX_j^2-x_jEX_j}{\sigma_j^2}], \hspace{5ex}
IF_{j,L_m}= \tilde r_m  \frac{x_j-EX_j}{\sigma_j^2}, \ \hspace{3ex} j=1,\ldots,p;
\end{equation}
$\sigma_j^2$ is the variance of $X_j, \ j=1,\ldots,p. $
\enp

%\pagebreak

{\em COMPARISON OF $L_m$-RESIDUALS FOR $F, \ F_{\epsilon, {\bf x},y}, \  F_{\epsilon, {\bf x}_{i,h},y},  F_{\epsilon, {\bf x},y+h }\  \ m=1,2$}

The next proposition confirms that for ${\bf x}$-remote case $({\bf x},y),$ the size of $L_1$ residual is  larger than the size of its $L_2$ residual
before $({\bf x},y)$ reaches $L_1$ {\em LBP}.
% {\em location} breakdown point.
\bep \label{p:residall}
For 
regression
 model (\ref{eq:mulregmod}) with  (${\cal A}1$)-(${\cal A}3$), perturbation (\ref{eq:xperturb}) and 
$r_1({\bf x},y) \neq 0:$\\
%and notation (\ref{eq:coef})- (\ref{eq:IFnotation}),
$a)$ For $\epsilon$ small: \\
$a_1)$  The difference of $({\bf x},y)$-residuals at $F_{\epsilon, {\bf x},y}$ and $F$ is:
\begin{equation}
\label{eq:residxadded}
 r_{m,{\bf x}}({\bf x},y) -r_m({\bf x},y) \approx -\epsilon [IF_{0,L_m}+\sum_{j=1}^px_jIF_{j,L_m}]
%[IF_{0,L_m}+x_1IF_{1,L_m}+\ldots+x_pIF_{p,L_m}]
=-\epsilon \tilde r_m({\bf x},y) [1+\sum_{j=1}^p\frac{(x_j-EX_j)^2}{\sigma_j^2}]; 
\end{equation}
 $r_{m,{\bf x}}({\bf x},y)$ and $r_{m}({\bf x},y)$ have the same sign and $|r_{m,{\bf x}}({\bf x},y)|<|r_{m}({\bf x},y)|,
 \ m=1, 2.$\\
$a_2)$ The ratio:
\begin{equation}
\label{eq:residxaddedL2L1ratio}
\frac{ r_{2,{\bf x}}({\bf x},y) -r_2({\bf x},y)}{ r_{1,{\bf x}}({\bf x},y) -r_1({\bf x},y)} \approx 2 \tilde f_{Y|{\bf X}}\frac{r_2({\bf x},y)}{sign[r_1({\bf x},y)] };
\end{equation}
$\tilde f_{Y|{\bf X}}$ is positive constant (\ref{eq:ftildmult}).\\
b) For $\epsilon$ and $|h|$ both small:\\
$b_1)$ The difference of $({\bf x},y)$-residuals at $F_{\epsilon, {\bf x},y}$ and $F_{\epsilon, {\bf x}_{i,h},y}$ is:
% $$ r_{m,{\bf x}_{i,h}}({\bf x}_{i,h},y)- r_{m,{\bf x}}({\bf x},y)+\beta_{i,L_m} h$$
 \begin{equation}
\label{eq:residxh}
%r_{m,{\bf x}_{i,h}}({\bf x},y)- r_{m,{\bf x}}({\bf x},y)+\beta_{i,L_m} h 
 r_{m,{\bf x}_{i,h}}({\bf x}_{i,h},y)- r_{m,{\bf x}}({\bf x},y)+\beta_{i,L_m} h
\approx -\epsilon h[IF_{i,L_m}+IF'_{x_i,0,L_m}+\sum_{j=1}^p x_j  IF'_{x_i,j, L_m}]
%[IF_{i,L_m}+ IF'_{x_i,0, L_m}+ x_1 IF'_{x_i,1, L_m}+x_2IF'_{x_i,2,  L_m}+ \ldots+ x_p IF'_{x_i,p, L_m}]
-\epsilon h^2  IF'_{ x_i, i, L_m}.
\end{equation}
%$$IF_{i,L_m}+IF'_{x_i,0,L_m}+\sum_{j=1}^p x_j  IF'_{x_i,j, L_m}$$
Thus,
\begin{equation}
\label{eq:r1dif}
 r_{1,{\bf x}_{i,h}}({\bf x}_{i,h},y)- r_{1,{\bf x}}({\bf x},y) \approx -\epsilon h \frac{sign[r_1({\bf x},y)]}{ \tilde f_{Y|{\bf X}}}\frac{x_i-EX_i}{\sigma_i^2} 
-\beta_{i,L_1} h -\epsilon h^2  IF'_{ x_i, i, L_1},
\end{equation}
%If ${\bf x}_{i,h}$ is more remote than ${\bf x},$ then $r_{1,{\bf x}_{i,h}}({\bf x}_{i,h},y), \  r_{1,{\bf x}}({\bf x},y)$ have both
%the same sign and $|r_{1,{\bf x}_{i,h}}({\bf x}_{i,h},y)|<| r_{1,{\bf x}}({\bf x},y)|.$ 
\begin{equation}
\label{eq:r2dif}
r_{2,{\bf x}_{i,h}}({\bf x}_{i,h},y)- r_{2,{\bf x}}({\bf x},y) \approx -\epsilon h \{ 2\frac{r_2(x_i-EX_i)}{\sigma_i^2}-\beta_{i,L_2} [1+\sum_{j=1}^p \frac{(x_j-EX_j)^2}{\sigma_j^2}]\}
 -\beta_{i,L_2} h -\epsilon h^2  IF'_{ x_i, i, L_2}.
\end{equation}
$b_2)$ If, in addition, $|x_i|$ is large,
\begin{equation}
\label{eq:r1difxrem}
 r_{1,{\bf x}_{i,h}}({\bf x}_{i,h},y)- r_{1,{\bf x}}({\bf x},y) \approx -\epsilon h \frac{sign[r_1({\bf x},y)]}{ \tilde f_{Y|{\bf X}}}\frac{x_i-EX_i}{\sigma_i^2},
%-\beta_{i,L_1} h -\epsilon h^2  IF'_{ x_i, i, L_1},
\end{equation}
 \begin{equation}
\label{eq:r2difxrem}
 r_{2,{\bf x}_{i,h}}({\bf x}_{i,h},y)- r_{2,{\bf x}}({\bf x},y) \approx \epsilon h  \cdot 3 \beta_{i,L_2}\frac{(x_i-EX_i)^2}{\sigma_i^2}, 
\end{equation}
%If ${\bf x}_{i,h}$ is more remote than ${\bf x},$ then $r_{m,{\bf x}_{i,h}}({\bf x}_{i,h},y), \  r_{m,{\bf x}}({\bf x},y)$ have %boththe same sign and $|r_{m,{\bf x}_{i,h}}({\bf x}_{i,h},y)|<| r_{m,{\bf x}}({\bf x},y)|, \ m=1,2.$ \\
%ii)
\begin{equation}
\label{eq:rratiodifxrem}
%0<
\frac{| r_{2,{\bf x}_{i,h}}({\bf x}_{i,h},y)- r_{2,{\bf x}}({\bf x},y)|}
{| r_{1,{\bf x}_{i,h}}({\bf x}_{i,h},y)- r_{1,{\bf x}}({\bf x},y)|}
\approx 3 |\beta_{i,L_2}| \cdot \tilde f_{Y|{\bf X}} \cdot
%\frac{
|x_i-EX_i|
%}{sign[r_1({\bf x},y)]}.
\end{equation}
{\em c)}  For $\epsilon$ and  $|h|$ both small, the difference of $({\bf x},y)$-residuals at $F_{\epsilon, {\bf x},y+h}$ and 
$F_{\epsilon, {\bf x},y}$ is:
\begin{equation}
\label{eq:r2dify}
r_{2,{\bf x},y+h}({\bf x},y+h)- r_{2, {\bf x},y} ({\bf x},y) \approx h-\epsilon h [1+\sum_{j=1}^p\frac{(x_j-EX_j)^2}{\sigma_j^2}].
\end{equation}
%{eq:sumIFderL2y}

 \enp

%{\em ${\bf x}$-INFLUENCE 
{\em INFLUENCE
%OF {\bf x}-PERTURBATIONS 
ON THE DERIVATIVES OF REGRESSION RESIDUALS-RINFIN}

Influence is determined from the distance  of residuals' 
%${\bf x}$-
derivatives at $({\bf x},y)$
% of residuals 
for model $F$
and gross-error model  $F_{\epsilon,{\bf x},y}.$
% \  F_{\epsilon,{\bf x}_{i,h},y}.$
The larger the distance is, the larger the influence of  $({\bf x},y)$ is.

{\em ${\bf x}$-Influence on Residuals}

For $({\bf x}_{i,h},y)$ and $({\bf x},y)$ both under model $F,$
\begin{equation}
\label{eq:rderFx}
\frac{r_m({\bf x}_{i,h},y)-r_m({\bf x},y)}{h}=-\beta_{i,L_m}, \hspace{4ex} i=1,\ldots,p, \ m=1,2.
\end{equation}

From the results for  gross-error models $F_{\epsilon,{\bf x},y}, \  F_{\epsilon,{\bf x}_{i,h},y},$ 
%it follows:
the  difference of residuals derivatives is obtained.
\bep
\label{p:index} For models $F,$ 
%and the gross-error models 
$F_{\epsilon,{\bf x},y}, \ F_{\epsilon,{\bf x}_{i,h},y}$  and $L_m$ regression it holds
\begin{equation}
\label{eq:index}
\lim_{h \rightarrow 0}  \frac{r_{m,{\bf x}_{i,h}}-r_{m,{\bf x}}}{h}+\beta_{i,L_m} =
 - \epsilon  [IF_{i,L_m}+IF'_{x_i,0,L_m}+\sum_{j=1}^p x_j  IF'_{x_i,j, L_m}], \hspace{2ex} i=1,\ldots,p,  \ m=1,2.
\end{equation}
\enp

%\ber
%When $p>1,$ coordinates other than the $i$-th are involved in $\sum_{j=1}^p x_j  IF'_{x_i,j, L_m}$ in (\ref{eq:index}),  %motivating the use of two  indices.
%\enr

From (\ref{eq:rderFx}) and (\ref{eq:index}),  the right side of the latter measures influence of ${\bf x}$'s  $i$-th coordinate
in the residual's derivative and provides the  motivation for defining influence. When $p>1,$ coordinates other than the $i$-th are involved in $\sum_{j=1}^p x_j  IF'_{x_i,j, L_m}$ in (\ref{eq:index}),  motivating the use of two influence indices.

\bef
\label{D:definfluence}
For gross-error model $F_{\epsilon, {\bf x}, y},$\\
a) the influence of ${\bf x}$'s $i$-th coordinate in the $L_m$-residual is
\begin{equation}
\label{D:ithinfluence}
\epsilon 
 \cdot |IF_{i,L_m}({\bf x},y)+IF'_{x_i,0,L_m}({\bf x},y)+\sum_{j=1}^p x_j  IF'_{x_i,j, L_m}({\bf x},y)|, \hspace{3ex}  m=1,2,
\end{equation}
b) the influence of ${\bf x}$  in the $L_m$-residual is
\begin{equation}
\label{D:iinfluence}
\epsilon \cdot  \sum_{i=1}^p |IF_{i,L_m}({\bf x},y)+IF'_{x_i,0,L_m}({\bf x},y)+\sum_{j=1}^p x_j  IF'_{x_i,j, L_m}({\bf x},y)|, \hspace{3ex}  m=1,2.
\end{equation}
\enf

Influences for models $F_{\epsilon_1,{\bf x}_1, y_1}, \  F_{\epsilon_2,{\bf x}_2, y_2}$ can be compared.

%\bef
%\label{d:i-thresinfluence}
%Case $({\bf x}_1,y_1)$ with weight $\epsilon_1$ is more influential 
% in its $i$-th coordinate for $L_m$-residuals 
%than case $({\bf x}_2,y_2)$ with weight $\epsilon_2$   if 
%$$
%\epsilon_2 \cdot |IF_{i,L_m}({\bf x}_2,y_2)+IF'_{x_i,0,L_m}({\bf x}_2,y_2)+\sum_{j=1}^p x_{2,j}  IF'_{x_i,j, L_m}({\bf x}_2,y_2)|
%$$
%\begin{equation}
%%\epsilon_2 \cdot |IF_{i,L_m}({\bf x}_2,y_2)+IF'_{x_i,0,L_m}({\bf x}_2,y_2)+\sum_{j=1}^p x_{2,j}  IF'_{x_i,j, L_m}({\bf x}_2,y_2)
%\le \epsilon_1 \cdot |IF_{i,L_m}({\bf x}_1,y_1)+IF'_{x_i,0,L_m}({\bf x}_1,y_1)+\sum_{j=1}^p x_{1,j}  IF'_{x_i,j, L_m}({\bf x}_1,y_1)|.
%\end{equation}
%\enf
 
\bef
\label{d:resinfluence}
Case $({\bf x}_1,y_1)$ with weight $\epsilon_1$ is more influential 
% in its $i$-th coordinate
for $L_m$-residuals than case $({\bf x}_2,y_2)$ with weight $\epsilon_2,$  $m=1,2,$   if 
$$
\epsilon_2 \cdot \sum_{i=1}^p |IF_{i,L_m}({\bf x}_2,y_2)+IF'_{x_i,0,L_m}({\bf x}_2,y_2)+\sum_{j=1}^p x_{2,j}  IF'_{x_i,j, L_m}({\bf x}_2,y_2)|
$$
\begin{equation}
%\epsilon_2 \cdot |IF_{i,L_m}({\bf x}_2,y_2)+IF'_{x_i,0,L_m}({\bf x}_2,y_2)+\sum_{j=1}^p x_j  IF'_{x_i,j, L_m}({\bf x}_2,y_2)
\le \epsilon_1 \cdot \sum_{i=1}^p |IF_{i,L_m}({\bf x}_1,y_1)+IF'_{x_i,0,L_m}({\bf x}_1,y_1)+\sum_{j=1}^p x_{1, j}  IF'_{x_i,j, L_m}({\bf x}_1,y_1)|.
\end{equation}
\enf

{\em The $L_2$-Residual Influence Index (RINFIN):} For gross-error model $F_{\epsilon, {\bf x}, y},$ 
 (\ref{D:iinfluence}) for $m=2$ becomes from (\ref{eq:sumIFderL2}) in the Appendix,
%For $({\bf x},y)$ with weight $\epsilon$ its $L_2$-influence index, obtained   via  (\ref{eq:sumIFderL2}), is the sum
%of influences for all coordinates:
\begin{equation}
\label{eq:globalindexL2}
{\mbox {\em RINFIN}}({\bf x},y;\epsilon,L_2)=\epsilon \cdot \sum_{i=1}^p \{|  \ 2\frac{r_2({\bf x},y)(x_i-EX_i)}{\sigma_i^2}-\beta_{i,L_2} [1+\sum_{j=1}^p \frac{(x_j-EX_j)^2}{\sigma_j^2}]\ | \}.
\end{equation}

{\em Abuse of notation:} Using {\em RINFIN}$({\bf x},y)$ instead of   {\em RINFIN}$({\bf x},y;\epsilon,L_2).$

Remote ${\bf x}$'s   have large {\em RINFIN}$({\bf x},y;\epsilon,L_2).$

\bep
\label{p:rinfinlimit} 
\begin{equation}
\label{eq:rinfinlimit}
\lim_{|x_i| \rightarrow \infty}  \mbox{RINFIN}({\bf x}, y;\epsilon, L_2)=\infty.
\end{equation}
\enp

\ber (RINFIN$^*)$
\label{r:rinfinstar} 
%When $p=1,$ the coefficient of $\beta_{i,L_2}$ in (\ref{eq:globalindexL2}) depends only on  ${\bf x}$'s $i$-th component.
%This is not the case when $p>1$  and
 To measure strictly the influence of ${\bf x}$'s $i$-th component, which is dominant
when $x_i$ is remote (see (\ref{eq:sumIFderL2xrem})),   
 use also:
\begin{equation}
\label{eq:globalindexL2star}
RINFIN^*({\bf x},y;\epsilon,L_2)=\epsilon \cdot \sum_{i=1}^p \{|  \ 2\frac{r_2({\bf x},y)(x_i-EX_i)}{\sigma_i^2}-\beta_{i,L_2} [1+ \frac{(x_i-EX_i)^2}{\sigma_i^2}]\ | \}.
\end{equation}
 Note that  RINFIN  dominates   $RINFIN^*$ in the simulations with the normal model  in section 4. However  $RINFIN^*$ allows to identify the  bad leverage cases in Hawkins-Bradu-Kass data. 
\enr

%Since $F$ is unknown, use (\ref{eq:globalindexL2}) for  data with $n$ cases  in two different ways:\\
%{\em a)} assume each case comes from $F$ and  compute its index (\ref{eq:globalindexL2}) with weight $\epsilon=1/n$,\\
%{\em b)} in the plot of absolute $L_1$ regression residuals against the square length of covariates, identify
%  remote plot points from $k$ neighboring
%cases as group. Assume these cases follow probability $G,$ component in the  mixture
%$$(1-\epsilon)F+\epsilon G,$$ with the means of $F$ and $G$ at large distance and $G$'s variance small.
 %Replace these $k$ cases by $(\bar {\bf x}_k,\bar y_k),$ with coordinates their averages and weight $\epsilon=k/n.$ Index 
%(\ref{eq:globalindexL2}) is calculated for the data consisting by this new case and data's $(n-k)$ remaining cases.

{\em $y$-Influence on Residuals}

Since $IF'_{y,j,L_1}$ vanishes for every $j,$ influence index from $y$-derivatives of residuals is only presented for 
$L_2$-regression.
 
For $({\bf x},y+h)$ and $({\bf x},y)$ both  under model $F,$
\begin{equation}
\label{eq:rderFy}
\frac{r_2({\bf x},y+h)-r_2({\bf x},y)}{h}=1, \hspace{4ex} i=1,\ldots,p.
\end{equation}
\bep
\label{p:indexy} For models $F,$ 
$F_{\epsilon,{\bf x},y}, \ F_{\epsilon,{\bf x},y+h}$  and $L_2$ regression it holds
\begin{equation}
\label{eq:indexy}
\lim_{h \rightarrow 0}  \frac{r_{2,{\bf x}, y+h}({\bf x},y+h)-r_{2,{\bf x},y}({\bf x},y)}{h}-1 
%r_{2,{\bf x},y+h}({\bf x},y+h)- r_{2, {\bf x},y} ({\bf x},y) -1
 \approx - \epsilon  [1+\sum_{j=1}^p\frac{(x_j-EX_j)^2}{\sigma_j^2}].
 \end{equation}
\enp

%From (\ref{eq:r2dify}) it follows that 
%\begin{equation}
%\label{eq:r2dify}
%\lim_{h \rightarrow 0}  \frac{r_{2,{\bf x}, y+h}({\bf x},y+h)-r_{2,{\bf x},y}({\bf x},y)}{h}-1 
% \approx - \epsilon  [1+\sum_{j=1}^p\frac{(x_j-EX_j)^2}{\sigma_j^2}].
% \end{equation}
\ber
\label{r:yresinfl}
From (\ref{eq:indexy}), the $y$-influence index 
% on the residuals derivatives
 is 
\begin{equation}
\label{eq:globalindexyL2}
\sum_{j=1}^p\frac{(x_j-EX_j)^2}{\sigma_j^2};
\end{equation}
it is maximized for  cases in  the extremes of the ${\bf x}$-coordinates 
%which
and  can be visually 
implemented 
%for some data  
with the proposed plot when 
%the 
${\bf x}$-coordinates have all the same sign. 
\enr

%\section{\bf The Plots}
\section{\bf Applications: {\em RINFIN} and  Residuals' Plots}

{\em RINFIN'S PERFORMANCE IN SIMULATIONS}

Data  $({\bf X}, Y)$ follows  
%$p$-dimensional 
 linear model  (\ref{eq:mulregmod}) 
with 
% parameter values  
$\beta=(1.5,.5,0,1,0,0,1.5,0,0,0,1, 0, \ldots, 0), \\ p \ge 11,$ 
and 
%for model $F$ 
with ${\bf X}$   obtained from    $p$-dimensional normal distribution,
${\cal N}({\bf 0}, {\bf \Sigma}),$ with $\Sigma_{i,j}=.5^{|j-i|}, 1 \le i, \ j \le p,$   as in  Alfons, Croux  and Gelper (2013, p.11). When $p<11,$  $\beta$'s first  $p$ coordinates are used and ${\bf X}$ is  obtained as in the case   $ p \ge 11.$
Contaminated {\bf X}  with  probability $G$ consists of $p$ independent normal random variables with mean $\mu$ and standard deviation 1. Each sample has size $n=200$ and various values for $p$ and $\mu$ are used, $p <n=200.$
Errors $e$ for model  $F$ and for $G$ are  independent, normal random variables with mean zero and variance 1. 
The percentage
of contaminated cases  $\epsilon=.10.$ The number of simulated samples $N=100.$

The 20 contaminated cases are the first cases in the data. {\em RINFIN} and {\em RINFIN}$^*$  values are  calculated and the
case numbers  of  the top 20 values  are recorded and compared with contaminated cases 1-20 to calculate the misclassification proportion. The results appear in  Table 1.  {\em RINFIN} misclassification proportion  is only reported being uniformly better than {\em RINFIN}$^*$  for the data used.

\begin{table}[h!]
\centerline{\begin{tabular}{|l|c|c|c|c|c|c|c|}\hline
$p$ &  $\mu=0$  & $\mu=.5$ & $\mu=1$ & $\mu=1.5$ & $ \mu=2$&$ \mu=2.5$ & $\mu=3$ \\\hline
2  & 0.9095 &0.866   & 0.7895  & 0.642  & 0.524   & 0.3805  & 0.281   \\ \hline
5 &0.91  & 0.857  & 0.7065   & 0.5035   & 0.2935   & 0.155  & 0.0685   \\ \hline
 10 &  0.9305  & 0.849  & 0.621   & 0.361   &0.1365   & 0.0435  & 0.0085  \\\hline
15 & 0.93    &  0.8485  & 0.558   & 0.242  &  0.062   & 0.01    & 0  \\ \hline
20 & 0.9085   & 0.8425   & 0.496    & 0.167  & 0.0335   & 0.0015   &  0 \\ \hline
60 & 0.933    & 0.755   & 0.2025    &  0.0075   & 0 & 0 &0  \\\hline 
100 & 0.928   & 0.6815   & 0.1005   &  0.001 & 0  & 0& 0 \\ \hline
140  &  0.9305   &  0.6275   & 0.048  & 0 & 0  & 0 & 0 \\ \hline
190 & 0.9105    & 0.7025  &  0.1535 & 0.0155   & 0.001  &0  & 5e-04  \\ 
    \hline\end{tabular}}
\caption{Average misclassification proportion  of  RINFIN
in $N=100$ simulated samples of size $n$=200, of which 20 are contaminated leverage cases,
from $p$-dimensional, normal-mixtures data with distance of means in each coordinate $\mu$ and variance 1.
}\label{table1}
\end{table}
 
The simulations follow the spirit  in Khan, Van Aelst and Zamar (2007). When there are no contaminated cases, i.e. $\mu=0,$ the misclassification proportion 
is near  90\% for all $p$-values. For any fixed $p$-value,  the misclassification proportion decreases as $\mu$  increases. For any fixed $\mu$-value the misclassification proportion decreases as $p$ increases except for an anomaly  when  $p=190$ due to its proximity to the sample size. By increasing $n$ to 250 cases, this anomaly disappears.
The {\em blessing} of high dimensionality is observed; for a theoretical confirmation see Yatracos (2013, Section 8, Proposition 8.1).

%{\em Reading Plots}
{\em READING RESIDUALS'  PLOTS}
 
\quad The goal is to identify quickly cases that do not follow the unknown model $F$ of the data's majority,
% in the data,
in particular  bad leverage cases.``Naive'' plots of  {\em absolute} residuals for $L_1$ and  $L_2$ regression against the  sum of squares of the independent variables are used.

%To determine bad leverage   cases   
Look for:\\
%at first 
$\bullet$  {\bf (A)} 
remote 
neighboring plot-points creating visual gaps in the $L_1$-plot's residuals but smaller gaps in the $L_2$-plot;
these are bad leverage cases 
% before reaching 
far from $L_1$  {\em LBP}. For a given ${\bf x},$ the  gap is ``large'' when the ratio of absolute residuals from the upper and lower  gap's borders is larger or equal to two.\\
 $\bullet$ {\bf (B)}
 a group of plot-points with {\em neighboring}
  horizontal axis projections, distant   from 
% those in
%cases in 
the
bulk of the plot, 
 with the $L_1$-absolute residuals forming a vertical strip and at least one  of them near zero;
% which have  
%$L_1$-absolute residuals visually  near zero;  
% indicating the corresponding case is
these are bad leverage cases 
 near  $L_1$  {\em LBP}.\\
% that is part of a cluster.\\
% and
%$\bullet$ 
%absolute residuals 
%some of them having very small values,
% in the plot, 
%especially those which along with
%neighboring cases on the horizontal axis form a remote cluster.
$\bullet$  {\bf (C)} If no unusual leverage cases are identified when plotting against the ${\bf x}$'s square length,  plot the
absolute residuals against each explanatory variable and  check whether there are remote ${\bf x}$-coordinates
for which  {(\bf A)}, {\bf (B)} hold. \\
$\bullet$  {\bf (D)} Large absolute residuals, especially  at the extremes of the ${\bf x}$-values in the data, indicating bad leverage
or other  outlying cases.

%\pagebreak

%{\em Using RINFIN with Data}
{\em USING RINFIN WITH DATA}

 The data 
$$D_n=\{ ({\bf x}_1, y_1), \ldots , ({\bf x}_n, y_n)\},  \hspace{8ex}  D_{n,-m}= D_n-\{({\bf x}_m,y_m)\}.$$  

To calculate  sample {\em RINFIN}$({\bf x}_m,y_m)$ estimate the  parameters in (\ref{eq:globalindexL2}) and 
 use $\epsilon=1/n:$ 
%and make  replacements in  (\ref{eq:globalindexL2}):\
{\em a)}  Use $D_{n,-m}$ to  obtain $L_2$-estimates $\hat \beta_{L_2}$ 
%of  linear regression coefficients
% $\hat \beta_{L_2}$ from data $D_{n,-m} $  and 
and 
%calculate the residual 
%estimate 
$\hat r_2({\bf x}_m,y_m).$\\
%to replace $r_2({\bf x},y).$\\
% in  (\ref{eq:globalindexL2}). \\
{\em b)} Estimate  $EX_i$ and $\sigma_i^2,$   respectively,   by the
sample average and sample variance ${\bf x}$-data's  $i$-th coordinate
% of the ${\bf x}$-data  
in $D_{n,-m}, \ i=1,\ldots,p.$\\
{\em c)} Use $\hat \beta_{L_2}$'s  $i$-th coordinate and replace $x_i$ with ${\bf x}_m$'s  $i$-th coordinate, $i=1,\ldots,p.$
% Replace $\beta_{i,L_2}$ and $x_i,$ respectively, by  the $i$-th coordinates of   $\hat \beta_{L_2}$  and  of  ${\bf x}_m;  \ i=1,%\ldots,p.$

If a group  $G$ of $k$ remote ${\bf x}$-neighboring  cases exists, 
$$G=\{ ({\bf x}_{i_1}, y_{i_1}), \ldots , ({\bf x}_{i_k}, y_{i_k})\} \subset D_n,$$
$D_n$   may follow a gross-error model.
Let  $\bar g$ be   the average of the elements in $G$ and use, instead of $D_n,$ new data
% replacing it  in $D_n.$
%The data set is now
$$(D_n-G) \cup \{\bar g\}.$$
Calculate {\em RINFIN}-values following {\em a)}-{\em c)}.  For  {\em RINFIN}$(\bar g)$ use
$\epsilon=k/n;$ in the remaining $(n-k)$ cases  weights are $1/n.$

%To calculate sample $RINFIN(\bar g)$ repeat {\em a)}-{\em c)} 
%with $\bar g$ instead of ${\bf x}_m,y_m)$ and $D_n-G \cup \bar g$ instead of $D_n$ and weight $k/n$  for $\bar g, \ 1/n$ 
%the remaining $n-k$ cases in $D_n-G \cup \bar g.$ 

With  $J$ groups, $G_1, \ldots, G_ J,$ of remote ${\bf x}$-neighboring cases, $G_k \cap G_l =\emptyset, k \neq l,$ obtain 
averages $\bar g_1, \ldots, \bar g_J,$
and use  data set
$$(D_n-\cup_{j=1}^JG_j)\cup \{\bar g_1, \ldots, \bar g_J\}.$$
Proceed with {\em a)}-{\em c)}.  For {\em RINFIN}$(\bar g_j)$ use $\epsilon_j=k_j/n, k_j$ is the cardinality of $G_j, \ j=1,\ldots,J;$
in the remaining cases weights are   $1/n.$

%\pagebreak

%{\em Data Plots}
{\em DATA PLOTS} \&  {\em RINFIN VALUES}

In Figures 1 and 2, $L_1$ and $L_2$ plots of absolute regression residuals are presented for twelve, well known data sets; those
without reference 
%may be found 
are in Rousseeuw and Leroy (1987). Several methods fail to determine 
cases from
% the 
gross-error component(s). The  visual comparison of regression plots, {\em RINFIN} and  {\em RINFIN}$^*$  are informative providing indications.
%which do not follow the model of the majority in the data.
%Potential unusual

Six  data sets  present   large  visual gaps  in $L_1$-plots  (see {(\bf A})) and smaller gaps in  $L_2$-plots.
%a gap is large when the ratio of absolute residuals from the upper and lower  gap's borders is larger or equal to two.
In one of the former,  the Hawkins-Bradu-Kass data,  residuals' plots give the impression of more than one  gross-error component but {\em RINFIN}$^*$ is successful after grouping.
%only data's creators know which one belong to the model $F.$ 
The remaining  data sets  presenting smaller $L_1$-gaps, if any, are:  Hertzsprung-Russel, Hadi-Simonoff, Stack Loss, Coleman, Salinity and Modified Wood. The ultimate data set  is the most challenging because there are no immediate visual indications, but
{\em RINFIN} is successful without grouping as well as after grouping along with {\em RINFIN}$^*.$
%both {\em RINFIN} and {\em RINFIN}$^*$  after grouping.

% Leverage cases with large absolute residuals are identified visually  in  several of the $L_1$-plots and are confirmed as bad %leverage cases by comparison  with  the $L_2$-plots where the size of the residuals is  much smaller.

% taking also into consideration the units. 
%In the Wood data the x-values seem to be either
%not ``semi-wild'' or beyond a location break-down point;
%the $L_1$ and $L_2$ absolute residuals are comparable.
%Studentized $L_1$ residual plots could be used to determine   theory becomes available; other confirmatory methods may also be %used.
%The $L_2$-plots
%are used to complement
%and compare with the findings of the $L_1$-plots.

In  Telephone data (covariates' dimension $p=1,$ number of cases $n=24,$ all covariates positive),
% page 26
 observations 15-20 cause a large gap in the residuals of the $L_1$-plot and no gap in the $L_2$-plot; 
${\bf (D)}$ 
applies. These are indeed 
the outliers because of the change in the recording system used.

 In  Kootenay river data ($p=1, \ n=13,$ all covariates positive),
%page 63
 case  4 is remote and 
% is a``bad''  leverage case and
%R&R page 64
  causes a large gap in the $L_1$-residuals, unlike the $L_2$-residuals. 
%It satisfies also the visual criterion for the $y$-influence index in the plot because all its covariates are positive. 
Both ${\bf (A)}$  and ${\bf (D)}$ apply.
{\em RINFIN} values confirm the visual findings.

 \begin{center}
\begin{tabular}{|c|c|c|c|c|c|c|}\hline
\multicolumn{7}{|c|}
{\bf  DATA: Kootenay River (p=1, n=13)}\\
 \hline CASE & 4  & 7& 2&  12  &  6 & 1\\
\hline {\em RINFIN}  & 8.906 & 0.106 &  0.052 & 0.044  & 0.030  & 0.015         \\
 \hline
\end{tabular}\\
\end{center}

%KOOTENAY DATA
 
% [8,]    1 1.886624e-01    .189    .189/13= 0.01453846 --> 0.015
% [9,]    6 3.874619e-01    .387                  .387/13= 0.02976923 -->0.030
%[10,]   12 5.703622e-01    .570           .570/13=0.04384615 ---> 0.044
%[11,]    2 6.726194e-01   .673           .673/13= 0.05176923  --> 0.052
%[12,]    7 1.382088e+00     1.382      1.382/13= 0.1063077 -->0.106
%[13,]    4 1.157730e+02    115.773  115.773/13= 8.905615  --> 8.906

% far from the bulk of the data and has large absolute residuals with projections forming  a gap among all projections in both axes %only in the $L_1$-plot. 

In  Brain and Body data ($p=1, \ n=28, $ not all covariates positive),
% page 57
 cases  6, 16, 25 are remote, obtained from 3 dinosaurs each with small brain and heavy body,
and  cause  a large gap  in the $L_1$ residuals. 
%Their  residual values are all larger in the $L_1$ plot. 
${\bf (A)}$  applies. {\em RINFIN}  values confirm the visual findings.
Case 26 in $R$ library is case 25 in Rousseeuw and Leroy (1987).

  \begin{center}
\begin{tabular}{|c|c|c|c|c|c|c|}\hline
\multicolumn{7}{|c|}
{\bf  DATA: LogBrain and LogBody (p=1, n=28)}\\
 \hline CASE & 25 &  6 & 16 & 27 & 17 & 10 \\
\hline {\em RINFIN}  & 0.298 &  0.183 & 0.157  & 0.051 &  0.039 &  0.029    \\
 \hline
\end{tabular}\\
\end{center}

%Log BRAIN AND Log BODY
%[23,]   10 0.80608675    0.80608675/28= 0.02878881 = 0.029
%[24,]   17 1.09490735     1.09490735/28= 0.03910383 = 0.039
%[25,]   27 1.41420897    1.41420897 /28= 0.05050746 =0.051
%[26,]   16 4.39758176    4.39758176/28= 0.1570565  =  0.157
%[27,]    6 5.12134852     5.12134852/28=0.1829053 = 0.183
%[28,]   26 8.33861054      8.33861054 /28=0.2978075 = 0.298

In  Hertzsprung-Russel star data ($p=1, \  n=47,$ all covariates positive),
%page 27
 cases 11, 20, 30, 34 correspond to giant stars.  
%Since $p=1,$ from the plots
 These are remote, $x$-neighboring cases
%far away from the bulk of the plot, 
%a large number 
and  many of the  remaining cases have
either comparable or  larger absolute residuals. Absolute residuals of  cases 11, 20, 30, 34  
form a narrow strip in the $L_1$ plot and that of case 11
is  near zero, indicating its proximity
 to $L_1$  {\em LBP}. Thus, ${\bf (B)}$ applies. 
Barrodale (1968, p. 55, l.
-2- p.56, l. 2) observed a similar behavior in an example for cases
he called ``wild''. 
%Most  $L_2$ absolute residuals are larger  than the corresponding  $L_1$ absolute  residuals;  (\ref{eq:residxaddedL2L1ratio}) %provides the explanation since  most  are smaller than 1.
%Note that $L_1$ and $L_2$ absolute residuals are in their vast majority smaller than one.
%Barrodale (1968, p. 55, l. -2- p.56, l. 2) observed a similar behavior in an example for cases he called ``wild''. 
%Moreover, % there is  a high  percentage of cases
%with unusually large absolute residuals was also observed.
{\em RINFIN} values  after grouping  ${\bf x}$-neighboring cases,
%  provide some information. However by considering, due to $x$-proximity, group 
$G_1=\{11, 20, 30,34\}, \ G_2=\{7,14\}$ 
support that $\{11, 20, 30,34\}$ are bad leverage cases.
Note that after grouping, the  cases in the table  form an ``envelope'' in the $L_1$ and $L_2$ plots.

\begin{center}
\begin{tabular}{|c|c|c|c|c|c|}\hline
\multicolumn{6}{|c|}
{\bf DATA:  Hertzsprung-Russel stars ($p=1, \  n=47$) }\\
% ($n=47, \ p=1,$ all covariates positive)\\
\hline CASE &{\em RINFIN} & GROUP & {\em RINFIN}& GROUP & {\em RINFIN}\\
\hline 34 &  0 .545  & 11,20,30,34 &    26.555 &   11,20,30,34 &    39.654   \\
\hline 30 &  0.387 &    14 &         0.276        &          7,14  &      0.447         \\
\hline 20& 0 .272 &  36  &       0.131      &         17 &       0 .159      \\
\hline 14& 0 .198 & 4  &   0.131  &                  36    &  0.149          \\
\hline 7 &   0.191 &   2   &    0.131   &   4     &      0.143        \\
\hline 11 & 0 .162 &  17  &    0 .125   &    2   &     0.143    \\
\hline
\end{tabular}\\
\end{center}

%RESULTS WITH WEIGHTS 1/47 ALL

%[42,]   11  7.6560249   7.6560249 *1/47=0.1628941   =.162
%[43,]    7  8.9797406    8.9797406 *1/47=0.1910583       =.191
%[44,]   14  9.3128202   9.3128202*1/47= 0.1981451       =.198
%[45,]   20 12.7835202   12.7835202*1/47= 0.2719898     =.272
%[46,]   30 18.2111254   18.2111254*1/47= 0.3874708     =.387
%[47,]   34 25.6022586   25.6022586*1/47= 0.5447289     = .545

%RESULTS WITH ONE GROUP 11,20,30,34, n=47-4+1=44

%[39,]   16   5.8776191    5.8776191  *1/47=0.1250557  =0.125
%[40,]    2   6.1432681     6.1432681 *1/47=0.1307078   =.131
%[41,]    4   6.1432681     6.1432681 *1/47=0.1307078    =.131
%[42,]   32   6.1504195   6.1504195 *1/47= 0.13086        =0.131
%[43,]   13  12.9575201   12.9575201 *1/47= 0.2756919  =0.276
%[44,]   44 312.0229485     312.0229485 *4/47= 26.55514   =  26.555

  %CORRESPONDANCE  44 -> 11,20,30,34
%13 --> 14, 
% 32 --> 36,
% 4 --> 4, 
%2 --> 2, 
%16 -->17 

%TAKING OUT TWO GROUPS NOW 11,20,30,34 (AT THE END AGAIN) AND  2nd group 7, 14    n=44-2+1=43

%[38,]    2   6.7361744     6.7361744 *1/47= 0.1433229 =.143
%[39,]    4   6.7361744     6.7361744 *1/47= 0.1433229  =.143
%[40,]   30   7.0045965    7.0045965*1/47=0.149034   =.149
%[41,]   14   7.4691986     7.4691986*1/47= 0.1589191   =.159
%[42,]   42  10.5035390    10.5035390*2/47=0.4469591  =  .447
%[43,]   43 465.9383957   465.9383957 * 4/47=39.65433  =39.654
 
%CORRESPONDANCE  43 -> 11,20,30,34
%42 --> 7,   14
%14 --> 17
% 30 -->  36
%  4   ->4 
% 2-->2 
In  
Hawkins-Bradu-Kass (1984)  artificial data ($p=3, \ n=75,$ all covariates positive),
%in the $L_1$-plot 
cases 11-14  have the largest absolute  residual in the $L_1$-plot  and are the most distant from cases 15-75.
The plot shows  three separated, distant groups that could be attributed to two sources of gross errors.
% It is a case where plotting against residuals would be more informative.
Using $({\bf A}),$ cases 11-14 are bad leverage cases. 
 Using {\bf (B)},  cases 1-10  are  ``bad'' leverage cases  near  $L_1$ {\em LBP}.
%indicating  in the $L_1$ plot  area of  location breakdown case.
Using two groups, cases $1-10$ and $11-14,$ {\em RINFIN}$^*$ indicates the true ``bad'' leverage cases $1-10.$
Rousseeuw and van Zomeren (1990, Figure 5)   identify cases 1-10 in the plot of standardized LMS residuals  against robust distances.

\begin{center}
\begin{tabular}{|c|c|c|c|c|c|}\hline
\multicolumn{6}{|c|}
{\bf DATA: Hawkins-Bradu-Kass ($p=3, \  n=75$) }\\
\hline CASE &{\em RINFIN} & GROUP & {\em RINFIN} & GROUP & {\em RINFIN}$^*$\\
\hline  12 &  0.523      & 11-14  &    17.848  & 1-10 & 16.648\\
\hline 14 &  0.442      &  1-10  &       15.983 & 11-14 &  14.426    \\
\hline  11 &   0.356   & 43     &   0.028   &43 &    0.028      \\
\hline 13 &   0.351    & 68 &        0.022  & 68  &  0.019         \\
\hline  7 &   0.174     &47          &  0.021   &  47 &   0.018         \\
\hline  6&  0.156       &27   &      0.019   &    27  &   0.017   \\
\hline  3  &  0.136     & 52       &   0.018  &   54     &   0.015   \\
\hline  5 &   0.133     & 60       &   0.0170   &  52   &  0.015    \\
\hline
\end{tabular}\\
\end{center}

 In  Scottish Hill Races data ($p=2, \ n=35,$ all covariates positive), cases 7 and 18 both at data's extremes cause large gap in the $L_1$-residuals unlike the $L_2$-residuals. $({\bf A})$ and ${(\bf D})$ both apply.
% where observation 18 causes the gap but is not a remote case.
%and 33 seem unusual; observation 7 is an unusual leverage point. 
According to Atkinson (1986, p. 399) observation 33 is masked by points 7 and 18 and an error is reported for  case 18.
%However, the comparison of $L_1$ and $L_2$
%plots indicates that case 33 is masked by case 7 because in the $L2$ plot the cases 35 and 33 have similar
%absolute residuals, that of case 7 is reduced but that of case 18 is unaltered.   
%Atkinson (1986) reports an error in case 18.
%it should be 18 minutes and not 78 minutes
{\em RINFIN}-values identify  these cases and in addition case 11.
% that can be barely seen in the plot in the line that would connect 0 with case 18. 

  \begin{center}
\begin{tabular}{|c|c|c|c|c|c|c|c|}\hline
\multicolumn{8}{|c|}
{\bf  DATA: Scottish Hill Races  (p=2, n=35)}\\
 \hline CASE & 7 &  11 & 35 & 33 & 18 & 31 & 17 \\
\hline {\em RINFIN}  &  3.577 & 3.230 &  2.492 &  1.558 & 1.067  &  0.796 &  0.508   \\
 \hline
\end{tabular}\\
\end{center}

%SCOTTISH HILL RINFIN READY!
%         CASE       RAW                  RINFIN
%[29,]   17      17.773061             0.5078017  =0.508
%[30,]   31       27.846926             0.7956265  = 0.796
%[31,]   18       37.331515             1.0666147 = 1.067
%[32,]   33       54.513637            1.5575325   = 1.558
%[33,]   35       87.226579             2.4921880  = 2.492
%[34,]   11      113.060753            3.2303072  =3.230
%[35,]    7       125.190308            3.5768659   = 3.577

In  Hadi-Simonoff (1993) data  ($p=2, n=25,$ one  covariate negative -.116), remote cases 1-3 have the larger absolute
residuals in the $L_1$ plot and case 17  follows. In the $L_2$
plot case 17 has the largest absolute residual and cases 1-3  
have their absolute residuals reduced. ${\bf (A)}$ applies for cases 1-3.
%After deletion of case 4 the $L_1$-plot shows clearly the
%``unusual'' nature of cases 1-3 but the picture does not change in the $L_2$-plot. 
When group $G=\{1,2,3\}$ is used, {\em RINFIN}  values indicate these are bad leverage cases.
% confirm the visual findings in the $L_1$-plot  when
%when group $G=\{1,2,3\}$ is used .
Hadi and Simonoff (1993) identify the {\em true outliers}, cases 1-3,
 and report that ``clean''
cases 6, 11, 13, 17 and 24 have larger absolute Least Median of
squares residuals than  cases 1-3. Plots of
standardized absolute residuals for bounded influence as well as
$M$-estimates regression do not reveal the outliers as unusual cases.
Yatracos (2013) identifies cases 1-3 as remote cluster with a projection-pursuit cluster  index.
%obtained from variance decompositions of the data's one-dimensional projections.
%$RINFIN$ values confirm the visual findings when, due to ${\bf x}$-proximity group 
%$G=\{1,2,3\}$ is used.

%The index is not informative because the gross-error distribution of cases 1-3 is in the boundary of the distribution $F;$  the 2-%dimensional factor space plot indicates that  cases 13, 17 and 4 are near cases 1-3. 

\begin{center}
\begin{tabular}{|c|c|c|c|}\hline
\multicolumn{4}{|c|}
{\bf  DATA:  Hadi-Simonoff ($p=2, \  n=25$) }\\
% ($n=47, \ p=1,$ all covariates positive)\\
\hline CASE  & {\em RINFIN} & GROUP  &{\em RINFIN}\\
\hline  22 &   0.677   & 1,2,3  &  1.074         \\
\hline 4 &   0.572 &  4  &     0.645            \\
\hline  17 &  0.527  & 17     &   0.620           \\
\hline 12 &  0.565  & 22 &  0.607                \\
\hline  25 &  0.374  &12      &  0.464                  \\
\hline  1&  0.351 &13   &     0.399              \\
\hline  2&    0.346   & 25 &    0.328         \\
\hline  3&    0.340   & 24 &    0.298          \\
\hline
\end{tabular}\\
\end{center}

 In 
Education data  ($p=3, \ n=50,$ all covariates positive), 
%page 110
case 50 (Alaska)
% is clearly different from the remaining cases in the $L_1$-plot.
causes a large gap in the $L_1$-residuals, unlike the $L_2$-residuals. 
${\bf (A)}$ applies. {\em RINFIN} values confirm the visual findings.

  \begin{center}
\begin{tabular}{|c|c|c|c|c|c|c|}\hline
\multicolumn{7}{|c|}
{\bf  DATA: Education   $(p=3, n=50)$}\\
 \hline CASE & 50 &  33 & 7 & 44 & 29 & 5  \\
\hline {\em RINFIN} & 1.20  &0.482  & 0.474  & 0.407  &  0.334 & 0.301      \\
 \hline
\end{tabular}\\
\end{center}

%EDUCATION DATA p=3, n=50
%         CASE       RAW                EPSILON ADJUSTED
%[45,]    5       15.0602881       0.30120576=0.301
%[46,]   29      16.7100362       0.33420072=0.334
%[47,]   44       20.3372884      0.40674577=0.407
%[48,]    7       23.7085651       0.47417130=0.474
%[49,]   33      24.1113104       0.48222621=0.482
% [50,]   50     59.9822083        1.19964417=1.20

 In  Stackloss data ($p=3, n=21,$ all covariates positive) 
%page 76
cases 1, 3, 4, 21  form a large gap  in the $L_1$-plot 
%rather
 and the gap  is reduced in the $L_2$-plot. ${\bf (A)}$ applies.
$L_1$-plot indicates   cases 1 and 3 can form a group with case 2 which has
small absolute residual and is near  {\em LBP}.  {\bf (B)} applies for cases 1,2,3.
%4 and  21  are unusual
% can be considered as  bad leverage cases. {\em RINFIN} and 
{\em RINFIN}$^*$ values indicate  $1,2, 3, 21$ are bad leverage cases. 
Case 4 has {\em RINFIN} and {\em RINFIN}$^*$ values before grouping at the 10-th percentile and after grouping 
below the 40th percentile. 
%Flores (2015, p. 805, 806)  identifies in  Stack Loss  data cases 1,2,3,4, 21 as bad leverage cases  using leverage constants. 
Rousseeuw and van Zomeren (1990) and   Flores (2015)  identify 
%in  Stack Loss  data 
cases 1, 2, 3, 4, 21 as bad leverage cases using, respectively,  plots of standardized LMS residuals against robust distances and   leverage constants. 

\begin{center}
\begin{tabular}{|c|c|c|c|c|c|c|c|}\hline
\multicolumn{8}{|c|}
{\bf  DATA:  Stackloss  ($p=3, n=21$) }\\
\hline CASE  &{\em RINFIN}& CASE  &{\em RINFIN}$^*$ & GROUP & {\em RINFIN} & GROUP & {\em RINFIN}$^*$\\
\hline 17  & 1.696   & 2  &    0.885 &  1,2,3   & 1.664   & 1,2,3    & 1.779\\
\hline  2 & 1.527   & 12  &   0.428   & 17    & 1.481  &   21      &  0.565      \\
\hline   1 & 0.757   & 21     &  0.427  &   21    & 0.697  &   12   &   0.444    \\
\hline  15& 0.557   & 17  &  0.420    &   7    &     0.642    &    7    &  0.409          \\
\hline   12 &  0.524 & 15     & 0.380    & 15   & 0.535     &   17    & 0.368                \\
\hline  18 &  0.520 &  11 & 0.317      &  8       & 0.531     &   15   & 0.358            \\
\hline  7 &   0.519    &  7 &    0.315     & 12    & 0.528  &  11     &   0.308         \\
\hline  8 &  0.440    &  16&    0.264       &  18   &  0.455  &  8     &  0.301       \\
\hline
\end{tabular}\\
\end{center}

In  Coleman data ($p=5, \ n=20,$ not all covariates positive)
% page 79 
cases 3 and 18  cause a large gap in the $L_1$-plot and ${\bf (D)}$ applies.
% and is located in the  extreme of the ${\bf x}$-data since all 
%${\bf x}$ coordinates are positive.
%${\bf (D)} applies.
Case  18 has larger absolute residual than most of the remaining cases and lives at the ${\bf x}$-extremes. 
%seem the most unusual in both
%The gaps in the $L_2$ residual are  visually larger but it is due to the plotting units.
%the $L_1$ plot and the $L_2$ plot; 
Cases 3, 4, 9, 16 in the $L_1$-plot indicate a potential cluster near  {\em LBP}.
%Cases 1 and 11 have the highest {\em RINFIN} and {\em RINFIN}$^*$ values and  those  for cases 15, 18, 19, 3, 16 follow.
 According to Rousseeuw and Leroy (1987)  `` .... examining the Least Squares results, ... cases 3, 11 and 18
are furthest away from the linear model. ...  The robust regression spots schools 3, 17 and 18 as outliers ...''.
{\em RINFIN} and {\em RINFIN}$^*$  highest four values identify  cases 1,6,11, 18, 19.

  \begin{center}
\begin{tabular}{|c|c|c|c|c|c|c|c|c|}\hline
\multicolumn{9}{|c|}
{\bf  DATA: Coleman (p=5, n=20)}\\
% \hline CASE & 1 & 11  & 15& 19   & 18 & 16  &13 & 17 \\
%\hline {\em RINFIN} & 3.399  &3.270  & 2.290  & 2.273   & 1.761 & 1.639  & 1.548 &1.445    \\
\hline CASE &  6   & 11  &  1 & 19   &10  & 15  &  2 & 16 \\
\hline {\em RINFIN} & 4.125  & 3.943   & 3.834  & 3.627  & 2.644  & 2.548  & 2.439  & 1.984    \\
 \hline CASE &11  & 1   &19  &  18  & 3 &12  & 6 &16 \\
\hline {\em RINFIN}$^*$  &2.297  & 1.736  & 1.316  & 1.291 & 1.279  & 1.269 & 1.150 & 1.076    \\
 \hline
\end{tabular}\\
\end{center}
 
% COLEMAN DATA RINFIN RESULTS: CASE, RAW WITHOUT EPSILON, 
%AND TIMES BY   1/n, n= 20  

%[13,]   17         28.893508      1.4446754=1.445
%[14,]   13        30.966780       1.5483390=1.548
%[15,]   16        32.781006       1.6390503=1.639
%[16,]   18        35.229128        1.7614564=1.761
%[17,]   19       45.450526         2.2725263=2.273
%[18,]   15       45.808767        2.2904383=2.290
%[19,]   11        65.393495       3.2696748=3.270
%[20,]    1          67.980799      3.3990400=3.399

%COLEMAN DATA RINFIN RESULTS: CASE, RAW, TIMES BY 1/n=1/20

%[13,]   16        21.528874            21.528874/20=1.076444=1.076
%[14,]    6         22.992418             22.992418/20=1.149621=1.150
%[15,]   12        25.371262             25.371262/20=1.268563=1.269
%[16,]    3         25.580564            25.580564/20= 1.279028=1.279
%[17,]   18         25.812697           25.812697/20= 1.290635=1.291
%[18,]   19         26.318199             26.318199/20=1.31591=1.316
%[19,]    1           34.715691            34.715691/20= 1.735785=1.736
%[20,]   11          45.942227             45.942227/20=2.297111=2.297

In  Salinity data (Ruppert and Carroll, 1980,  $p=3, \ n=28,$ all covariates positive)
% page 82
  case 16 has the largest absolute residual,
 is ${\bf x}$-remote and the gap caused  in the $L_1$ plot is small. In the $L_2$ plot its absolute residual is reduced.
 Both $({\bf A})$ and $({\bf D})$ apply for case 16.
%The least median of squares identifies both cases 5 and 16 as outliers.
%unusual in the $L_1$ plot. 
{\em RINFIN}  values confirm the visual findings. In Carroll and Ruppert (1985) the analysis of the data shows that cases
3 and 16 are masking case 5.

  \begin{center}
\begin{tabular}{|c|c|c|c|c|c|c|}\hline
\multicolumn{7}{|c|}
{\bf  DATA: Salinity (p=3, n=28)}\\
 \hline CASE & 16 & 15  & 5 &  3 & 9 & 4 \\
\hline {\em RINFIN} & 2.216 &  0.418 & 0.327  & 0.307 & 0.293  &  0.288    \\
 \hline
\end{tabular}\\
\end{center}

%   SALINITY DATA RINFIN p=3, n=28

%             CASE              RAW            WITH EPSILON
%[23,]        4             8.055685          0.28770302=0.288
%[24,]         9             8.203677         0.29298845=0.293
%[25,]         3             8.587905          0.30671089=0.307
%[26,]          5            9.148155          0.32671980=0.327
%[27,]          15          11.714650        0.41838034=0.418
%[28,]          16          62.051657         2.21613061=2.216

In  modified Wood data ($p=5, n=20,$ all covariates positive)
% page 243
there is no visual gap in the $L_1$-plot. 
Since covariates are positive, cases 7, 19, 1, 4, 6, 8  live  at one ${\bf x}$-extreme of the data and 
present a pattern like that described in ${\bf (B)},$
with the residual of case 8 near 0. To determine the  neighboring cases, plot $L_2$ (or $L_1$) residuals
against each ${\bf x}$-coordinate.
In Figure 3 it is clear that cases 4, 6, 8, 19 are neighboring and remote in each coordinate.
Cases 4, 6, 8, 19 form also a strip in the last two $L_2$-plots  in Figure 3. ${\bf (B)}$ applies for
these cases in view of the $L_1$ plot in Figure 1. This is confirmed by  the four higher  {\em RINFIN} values  and
also by both {\em RINFIN} and  {\em RINFIN}$^*$ values when $\{4,6,8,19\}$ are considered as group.

\begin{center}
\begin{tabular}{|c|c|c|c|c|c|c|c|}\hline
\multicolumn{8}{|c|}
{\bf  DATA:  Modified Wood  ($p=5, n=20$) }\\
\hline CASE  & {\em RINFIN} & CASE  &{\em RINFIN}$^*$                 & GROUP & {\em RINFIN} & GROUP & {\em RINFIN}$^*$\\
\hline  19 &  1.579 &    11   &   0.871                                     &    4,6,8,19   &  34.729                 &    4,6,8,19& 29.583    \\
\hline    8  &1.532   & 12    & 0.508                                           &    11    & 1.710                                         & 11   &  1.366      \\
\hline   6   &   1.452 &  1    &  0.493                                        &   7   &  1.460                                       & 7     &  0.597      \\
\hline  4   & 1.3312  &  7    &  0.476                                       &    12   & 1.390                                                   &   12  &   0.556      \\
\hline    12    &   1.324   &14    &   0.448                                    &   10   & 1.084                                                   & 1   &   0.468       \\
\hline   11  & 1.161    &   19   &  0.442                                        &   16   &    0.785                                        & 14   &   0.389        \\
\hline    7 & 1.158        &    8   & 0.434                                        &   1    &    0.779                                          &  16   & 0.285           \\
\hline   10    &  1.075     &   4   & 0.386                                         &   17 &   0.738	                                        &   10   &  0.252         \\
\hline
\end{tabular}\\
\end{center}

\section{APPENDIX}

{\bf Proof of Lemma \ref{l:key}:}
 Equality  (\ref{eq:limit1})  is obtained by adding and
subtracting $T(F)$
 in the numerator of its left side
%(\ref{eq:locbreakdown})
%(\ref{eq:limit1})
 and by taking first the
limit with respect to $\epsilon. \hspace{3ex}  \Box$

{\bf Proof for Proposition \ref{p:Elenimatrix}:} {\em a)} Induction is used. \\
For $n=1$, the determinant is $A_1-a_1^2.$\\
For $n=2,$ the determinant is
$$ (A_1 A_2-a_1^2a_2^2) -a_1 \cdot (a_1A_2-a_1a_2^2)+a_2 \cdot (a_1^2a_2-A_1a_2)
=A_1A_2-a_1^2A_2 +a_1^2a_2^2-A_1a_2^2$$
$$=A_2(A_1-a_1^2)-a_2^2(A_1-a_1^2)=(A_1-a_1^2) (A_2-a_2^2).$$ 
Assume that (\ref{eq:ElenimatrixDet}) holds for $E_n.$ 
%It is shown that it also holds for $n=k+1.$
To show it holds for $E_{n+1}$ consider the matrix $E_{n+1}:$
%\pagebreak
$$
 E_{n+1}=
 \left  ( 
\begin{array}{ccccc}
1 & a_1 & a_2 \ldots  & a_n& a_{n+1}\\
a_1 & A_1 & a_1 a_2  \ldots & a_1 a_n & a_1 a_{n+1}\\
a_2 & a_2a_1 & A_2 \ldots & a_2a_n  & a_2a_{n+1}\\
\ldots \\
a_n & a_na_1 & a_na_2 \ldots  & A_n &  a_n a_{n+1}\\
a_{n+1} & a_{n+1}a_1 & a_{n+1}a_2 \ldots  & a_{n+1}a_n &  A_{n+1}.\\
\end {array} 
\right )
$$
$|E_{n+1}|$ is  obtained using  line $(n+1)$  and its cofactors $C_{n+1,1}, \ldots, C_{n+1,n+1}:$
\begin{equation}
\label{eq:detE1}
|E_{n+1}|= a_{n+1}C_{n+1,1}+ a_{n+1}a_1  C_{n+1,2}+\ldots+ a_{n+1}a_n C_{n+1,n}+A_{n+1} C_{n+1,n+1}.
\end{equation}
Observe that for $2 \le j \le n,$  cofactor $C_{n+1,j}$  is obtained from a matrix where the last column is a multiple
of its  first    column  by $a_{n+1},$ thus,
\begin{equation}
\label{eq:cofactorsalmostall}
C_{n+1,j}=0, \ j=2,\ldots, n.
\end{equation} 
For the matrix in cofactor $C_{n+1,1},$ observe that in its  last column   $a_{n+1}$    is common factor and if taken out of the determinant 
the remaining column is the vector generating $E_n,$  i.e. $\{1,a_1,\ldots,a_n\}.$ With $n-1$ successive interchanges to the left, this column becomes first
and $E_{n}$ appears.
%moved out of the determinant,  being common factor of each element in the  last column, 
%then $E_k$'s determinant appears after interchanging the altered   last column  with the columns to its left
 %in  $k-1$ successive interchanges.
Thus,
 \begin{equation}
\label{eq:cofactor1}
C_{n+1,1}=(-1)^{n+2}(-1)^{n-1} \cdot a_{n+1}|E_n|=-a_{n+1}|E_n|.
\end{equation} 
In  cofactor $C_{n+1, n+1},$ the determinant is that of $E_n,$
  \begin{equation}
\label{eq:cofactor1}
C_{n+1,n+1}=(-1)^{2(n+1)} |E_n|=|E_n|.
\end{equation} 
From  (\ref{eq:detE1})-(\ref{eq:cofactor1}) it follows that
$$|E_{n+1}|=-a_{n+1}^2|E_n|+A_{n+1}|E_n|=\Pi_{m=1}^{n+1}(A_m-a_m^2). $$

{\em b)} We now work with $E_n.$  For $ i>0, j>0, i \neq j,$  after deleting row $(j+1)$  the remaining of column $(j+1)$ in the cofactor  is a multiple of column 1, thus
$|C_{i+1,j+1}|$ vanishes.

For $C_{1,j+1},$ using column $j+1$ to calculate $E_n,$ it holds:
$$a_j C_{1,j+1}+A_j C_{j+1,j+1}=|E_n| \rightarrow  a_j C_{1,j+1}=-a_j^2 \Pi_{k \neq j} (A_k-a_k^2) \rightarrow  C_{1, j+1}=-a_j  \Pi_{k \neq j} (A_k-a_k^2).$$
For $C_{i+1,1}, i>0,$  after deletion of row  $(i+1)$ in $E_n$  the remaining of  column $(i+1)$in the cofactor's matrix  is multiple of $a_i$  and  the basic vector creating $E_{n, -i}.$ 
%except for $a_i.$ 
Column 1 of $E_n$  is also deleted and for column $(i+1)$ in the cofactor's matrix  to become first column $(i-1)$ exchanges of columns are needed. Thus,
$$C_{i+1,1}=(-1)^{i+2} \cdot  a_i \cdot (-1)^{i-1} \Pi_{k \neq i}(A_k-a_k^2)=-a_i \cdot \Pi_{k \neq i}(A_k-a_k^2).$$
For $C_{1,1}$ we express $|E_n|$ as sum of cofactors along the first row of $E_n,$
$$C_{1,1}+a_1 C_{1,2}+\ldots+a_n C_{1,n}=|E_n| $$
$$\rightarrow C_{1,1}=\Pi_{k=1}^n(A_k-a_k^2) +
a_1^2 \Pi_{k\neq 1}(A_k-a_k^2) +\ldots +a_n^2  \Pi_{k\neq n}(A_k-a_k^2).  \hspace{5ex}\Box  $$

{\bf Proof of Proposition \ref{p:IF}:} For system of equations (\ref{eq:syst1}), (\ref{eq:syst2})
and  matrix $E_p$ with $a_j=EX_j, \ A_j=EX_j^2, \ j=1,\ldots, p,$ from Proposition \ref{p:Elenimatrix} 
$$IF_{j,L_m}=\frac{C_{1,j+1}\tilde  r_m +C_{j+1,j+1}\tilde  r_m x_j}{|E_p|}=\tilde  r_m\frac{-EX_j \Pi_{k \neq j}\sigma_k^2+x_j \Pi_{k \neq j}\sigma_k^2}
{\Pi_{k=1}^p \sigma_k^2}=\tilde  r_m \frac{x_j-EX_j}{\sigma_j^2}, \ j=1,\ldots,p.$$

$$IF_{0,L_m}=\frac{C_{1,1} \tilde  r_m +\sum_{j=1}^p C_{1,j+1} \tilde  r_m x_j}{|E_p|}=\tilde  r_m \frac{\Pi_{k=1}^p \sigma_j^2+\sum_{j=1}^p
 (EX_j)^2\Pi_{k \neq j}\sigma_k^2 - \sum_{j=1}^p x_j EX_j \Pi_{k \neq j}\sigma_k^2} {\Pi_{k=1}^p\sigma_k^2}$$
$$=\tilde  r_m[1+\sum_{j=1}^p \frac{EX_j^2-\sigma_j^2-x_jEX_j}{\sigma_j^2}=\tilde  r_m [1-p+\sum_{j=1}^p \frac{EX_j^2-x_jEX_j}{\sigma_j^2}]. \hspace{5ex} \Box$$ 

\bel
\label{eq:sumderIF}
For the influence functions (\ref{eq:IFALL}) it holds: \\
a)  
\begin{equation}
\label{eq:sumIFLm}
 IF_{0,L_m}+\sum_{j=1}^px_jIF_{j,L_m}=\tilde r_m [1+\sum_{j=1}^p\frac{(x_j-EX_j)^2}{\sigma_j^2}], \ m=1,2,
\end{equation}
b) 
\begin{equation}
\label{eq:sumIFderL1}
IF_{i,L_1}+IF'_{x_i,0,L_1}+\sum_{j=1}^p x_j  IF'_{x_i,j, L_1}=\frac{sign[r_1({\bf x},y)]}{ \tilde f_{Y|{\bf X}}}\frac{x_i-EX_i}{\sigma_i^2},
\end{equation}
c)
\begin{equation}
\label{eq:sumIFderL2}
IF_{i,L_2}+IF'_{x_i,0,L_2}+\sum_{j=1}^p x_j  IF'_{x_i,j, L_2} =2\frac{r_2({\bf x},y)(x_i-EX_i)}{\sigma_i^2}-\beta_{i,L_2} [1+\sum_{j=1}^p \frac{(x_j-EX_j)^2}{\sigma_j^2}]
\end{equation}
\begin{equation}
\label{eq:sumIFderL2xrem}
 \approx  -3 \beta_{i,L_2}\frac{(x_i-EX_i)^2}{\sigma_i^2}, \mbox{ if } |x_i-EX_i| \mbox{ is very large,}
\end{equation}
d) 
\begin{equation}
\label{eq:sumIFderL2y}
IF'_{y,0,L_2}+\sum_{j=1}^p x_j  IF'_{y,j, L_2} = 1+\sum_{j=1}^p \frac{(x_j-EX_j)^2}{\sigma_j^2}.  
\end{equation}
  \enl

{\bf Proof of Lemma \ref{eq:sumderIF}:}
{\em a)}  From (\ref{eq:IFALL}),
$$ IF_{0,L_m}+\sum_{j=1}^px_jIF_{j,L_m}= \tilde r_m[1-p+\sum_{j=1}^p \frac{EX_j^2-x_jEX_j}{\sigma_j^2}] + \sum_{j=1}^p x_j \frac{\tilde r_m(x_j-EX_j)}{\sigma_j^2}$$
$$=\tilde r_m [1-p + \sum_{j=1}^p \frac{EX_j^2-2x_jEX_j+x_j^2}{\sigma_j^2}]=\tilde r_m[1+\sum_{j=1}^p\frac{(x_j-EX_j)^2}{\sigma_j^2}], \ m=1,2.$$
 {\em b)}  Proof is  provided for $i=1.$  If the residual of $({\bf x},y)$ does not vanish, since
 $$IF_{0, L_1}= \frac{sign[r_1({\bf x},y)]}{2 \tilde f_{Y|{\bf X}}} [1-p+\sum_{j=1}^p \frac{EX_j^2-x_jEX_j}{\sigma_j^2}], \hspace{5ex}   IF_{j,L_1}=\frac{sign[r_1({\bf x},y)]}{2 \tilde f_{Y|{\bf X}}}\frac{x_j-EX_j}{\sigma_j^2}, \ j=1,\ldots,p,$$
$$IF'_{x_1,0,L_1}=- \frac{sign[r_1({\bf x},y)]}{2 \tilde f_{Y|{\bf X}}}\frac{EX_1}{\sigma_1^2}$$
$$IF'_{x_1,1,L_1}=\frac{sign[r_1({\bf x},y)]}{2 \tilde f_{Y|{\bf X}} \sigma_1^2}, \hspace{5ex} IF'_{x_1,j,L_1}=0, \ j \neq 1. $$
Thus,
$$ IF_{1,L_1}+ IF'_{x_1,0, L_1}+ x_1 IF'_{x_1,1, L_1}+x_2IF'_{x_1, 2, L_1}+ \ldots+ x_p IF'_{ x_1,p,  L_1}$$
$$=\frac{sign[r_1({\bf x},y)]}{2 \tilde f_{Y|{\bf X}}}\frac{x_1-EX_1}{\sigma_1^2}- \frac{sign[r_1({\bf x},y)]}{2 \tilde f_{Y|{\bf X}}}\frac{EX_1}{\sigma_1^2}+x_1 \frac{sign[r_1({\bf x},y)]}{2 \tilde f_{Y|{\bf X}} \sigma_1^2}=\frac{sign[r_1({\bf x},y)]}{ \tilde f_{Y|{\bf X}}}\frac{x_1-EX_1}{\sigma_1^2} $$
{\em c)}  Proof is  provided for $i=1.$ Since 
$$IF_{0, L_2}=  r_2[1-p+\sum_{j=1}^p \frac{EX_j^2-x_jEX_j}{\sigma_j^2}], \hspace{5ex}   IF_{j,L_2}= r_2 \frac{x_j-EX_j}{\sigma_j^2}, \ j=1,\ldots,p,$$
$$IF'_{x_1,0,L_2}=-\beta_{1,L_2} [1-p+\sum_{j=1}^p \frac{EX_j^2-x_jEX_j}{\sigma_j^2}]-r_2\frac{EX_1}{\sigma_1^2}$$
 $$IF'_{x_1,1,L_2}=-\beta_{1,L_2} \frac{x_1-EX_1}{\sigma_1^2}+\frac{r_2}{\sigma_1^2} \rightarrow 
x_1IF'_{x_1,1,L_2}=-\beta_{1,L_2} \frac{x_1^2-x_1EX_1}{\sigma_1^2}+r_2\frac{ x_1}{\sigma_1^2}$$
$$IF'_{x_1,j,L_2}=-\beta_{1,L_2} \frac{x_j-EX_j}{\sigma_j^2} \rightarrow x_j IF'_{x_1,j,L_2}=-\beta_{1,L_2} 
 \frac{x_j^2-x_jEX_j}{\sigma_j^2}, \hspace{3ex}  \ j \neq 1.$$
Thus,
$$IF_{1,L_2}+ IF'_{x_1,0, L_2}+ x_1 IF'_{x_1,1, L_2}+x_2IF'_{ x_1,2, L_2}+ \ldots+ x_p IF'_{x_1,p, L_2}$$
$$=2\frac{r_2(x_1-EX_1)}{\sigma_1^2}-\beta_{1,L_2}[1-p+\sum_{j=1}^p \frac{x_j^2-2x_jEX_j+EX_j^2}{\sigma_j^2}]$$
$$=2\frac{r_2(x_1-EX_1)}{\sigma_1^2}-\beta_{1,L_2} [1+\sum_{j=1}^p \frac{(x_j-EX_j)^2}{\sigma_j^2}]. $$
Since
$$r_2(x_1-EX_1)=y(x_1-EX_1)-\beta_{1,L_2}x_1(x_1-EX_1)-(x_1-EX_1)\sum_{j=2}^p \beta_{j,L_2}x_j$$
$$=y(x_1-EX_1)-\beta_{1,L_2}(x_1-EX_1)^2-\beta_{1,L_2}EX_1(x_1-EX_1)-(x_1-EX_1)\sum_{j=2}^p \beta_{j,L_2}x_j,$$
if $|x_1-EX_1|$ is very large dominating all the other terms, then
$$IF_{1,L_2}+ IF'_{x_1,0, L_2}+ x_1 IF'_{x_1,1, L_2}+x_2IF'_{2, x_1, L_2}+ \ldots+ x_p IF'_{p, x_1, L_2} \approx -3 \beta_{1,L_2}\frac{(x_1-EX_1)^2}{\sigma_1^2}. $$
{\em d)}   From (\ref{eq:IFALL}),
%$$IF_{0,L_m}=\tilde r_m [1-p+\sum_{j=1}^p \frac{EX_j^2-x_jEX_j}{\sigma_j^2}], \hspace{5ex}
%IF_{j,L_m}= \tilde r_m  \frac{x_j-EX_j}{\sigma_j^2}, \ \hspace{3ex} j=1,\ldots,p;$$
$$IF'_{y,0,L_2}=1-p+\sum_{j=1}^p \frac{EX_j^2-x_jEX_j}{\sigma_j^2}, \hspace{5ex} IF'_{y,j,L_2}= \  \frac{x_j-EX_j}{\sigma_j^2}, \ j=1,\ldots,p.$$
Thus,
$$IF'_{y,0,L_2}+\sum_{j=1}^p x_j IF'_{y,j,L_2}=1-p+\sum_{j=1}^p\frac{EX_j^2-x_jEX_j+x_j^2-x_jEX_j}{\sigma_j^2}=1+\sum_{j=1}^p\frac{(x_j-EX_j)^2}{\sigma_j^2}. \hspace{3ex} \Box$$

A Lemma used repeatedly  to calculate residuals'differences  is due to (\ref{eq:firstorder}), (\ref{eq:addedvalue}).

\bel 
\label{l:coefapprox}
For regression model (\ref{eq:mulregmod}) with assumptions (${\cal A}1$), (${\cal A}3$), perturbation (\ref{eq:xperturb}),
$r_1({\bf x},y) \neq 0,$ and
$\epsilon, |h|$ both small: 
%and notation (\ref{eq:coef})- (\ref{eq:IFnotation}),
\begin{equation}
\label{eq:coefapprox}
\beta_{j, L_m,{\bf x}} \approx \beta_{j,L_m}+\epsilon IF_{j,L_m}, \hspace{8ex}
\beta_{j, L_m,{\bf x}_{i,h}}\approx \beta_{j, L_m,{\bf x}}+ \epsilon h IF'_{x_i,j,L_m}.
\end{equation}
\enl

{\bf Proof of Lemma \ref {l:coefapprox}:}  Use approximations (\ref{eq:firstorder}), (\ref{eq:addedvalue}). $\hspace{5ex}\Box$

{\bf Proof of Proposition \ref{p:residall}:} {\em a)} For $a_1),$  from Lemma \ref{l:coefapprox},
 $$ r_{m,{\bf x}}=y-\beta_{0,L_m,{\bf x}}-\beta_{1,L_m,{\bf x}}x_1-\ldots -\beta_{p,L_m,{\bf x}}x_p$$
$$\approx y-(\beta_{0,L_m}+\epsilon IF_{0,L_m})-(\beta_{1,L_m}+\epsilon IF_{1,L_m})x_1-\ldots - 
(\beta_{p,L_m}+\epsilon IF_{p,L_m})x_p$$
$$=r_m-\epsilon(IF_{0,L_m}+x_1 IF_{1,L_m}+\ldots+ x_p IF_{p,L_m}).$$
(\ref{eq:residxadded}) follows from (\ref{eq:sumIFLm}). Since $\tilde r_m({\bf x},y)$ has the same sign with $ r_m({\bf x},y),$ for $\epsilon$ small
$r_{m,{\bf x}}({\bf x},y)$ will also have the same sign and reduced size because  $-\epsilon \tilde  r_m({\bf x},y)$ has opposite sign from  $ r_m({\bf x},y).$\\
 For $a_2),$ (\ref{eq:residxaddedL2L1ratio}) follows from (\ref{eq:rtilde}).\\
{\em b)}   Provided for $i=1$ using  Lemma \ref{l:coefapprox}:
 $$r_{m,{\bf x}_{1,h}}=y-\beta_{0,L_m,{\bf x}_{1,h}}-\beta_{1,L_m,{\bf x}_{1,h}}(x_1+h)-\ldots - \beta_{p,L_m,{\bf x}_{1,h}}x_p$$
$$\approx y-[\beta_{0,L_m,{\bf x}} + \epsilon h IF'_{x_1,0, L_m}]-[\beta_{1,L_m,{\bf x}} + \epsilon h IF'_{ x_1, 1, L_m}](x_1+h)
-\ldots - [\beta_{p,L_m,{\bf x}} + \epsilon h IF'_{x_1, p, L_m}]x_p$$
$$=r_{m,{\bf x}}-\beta_{1,L_m,{\bf x}} h-\epsilon h[IF'_{ x_1,0, L_m}+x_1 IF'_{ x_1,1,  L_m}+x_2IF'_{x_1,2,  L_m}+ \ldots+ x_p IF'_{x_1,p, L_m}]-\epsilon h^2  IF'_{x_1,1,  L_m}$$
%$$=r_{m,x_1}-\beta_{1,L_m,x_1} h -\epsilon h[IF'_{0, x_1, L_m}+x_1 IF'_{1, x_1, L_m}+x_2IF'_{2, x_1, L_m}+ \ldots+ x_p IF'_{p, x_1, L_m}]-\epsilon h^2  IF'_{x_1,1, L_m}$$
$$=r_{m,{\bf x}}-\beta_{1,L_m} h - \epsilon h[IF_{i,L_m}+ IF'_{x_1,0, L_m}+ x_1 IF'_{x_1,1, L_m}+x_2IF'_{2, x_1, L_m}+ \ldots+ x_p IF'_{p, {\bf x}, L_m}]-\epsilon h^2  IF'_{ x_1, 1, L_m}.$$
(\ref{eq:r1dif}), (\ref{eq:r2dif}) follow from (\ref{eq:sumIFderL1}), (\ref{eq:sumIFderL2}).
%(\ref{eq:sumIFderL2xrem})

For $b_2),$ if $|x_i|$ is large and $|h|$ is small, $\beta_{i,L_m}h$ and $\epsilon h^2 IF'_{x_i,i,L_m}$ are of smaller order than the remaining terms and (\ref{eq:r1difxrem}) follows, in addition,  (\ref{eq:sumIFderL2xrem}) implies (\ref{eq:r2difxrem}) and 
(\ref{eq:rratiodifxrem}) follows also.
 $$ c) \hspace{5ex} r_{2,{\bf x}, y+h}({\bf x},y+h)=y+h-\beta_{0,L_2,{\bf x}, y+h}-\beta_{1,L_2,{\bf x},y+h}x_1-\ldots - \beta_{p,L_2,{\bf x },y+h}x_p$$
$$\approx y+h-[\beta_{0,L_2,{\bf x},y} + \epsilon h IF'_{y,0, L_2}]-[\beta_{1,L_2,{\bf x},y} + \epsilon h IF'_{ y, 1, L_2}]x_1
-\ldots - [\beta_{p,L_2,{\bf x},y} + \epsilon h IF'_{y, p, L_2}]x_p$$
$$=r_{2,{\bf x},y}({\bf x},y) + h-\epsilon h [IF'_{y,0, L_2}+\sum_{j=1}^p x_j IF'_{ y, j, L_2}]=r_{2,{\bf x},y}({\bf x},y) + h-\epsilon h [1+\sum_{j=1}^p\frac{(x_j-EX_j)^2}{\sigma_j^2}],$$
with the last equality obtained from (\ref{eq:sumIFderL2y}). $\hspace{5ex} \Box$

{\bf Proof of Proposition \ref{p:index}:} Follows from (\ref {eq:residxh}) dividing both its sides by $h$ and taking the limit with $h$ converging to zero. $\hspace{2ex} \Box$

{\bf Proof of Proposition \ref{p:rinfinlimit}:}  
$$\lim_{|x_i| \rightarrow \infty} \mbox{{\em RINFIN}}({\bf x},y;\epsilon,L_2) \ge \epsilon \cdot  \lim_{|x_i| \rightarrow \infty}|  \ 2\frac{r_2({\bf x},y)(x_i-EX_i)}{\sigma_i^2}-\beta_{i,L_2} [1+\sum_{j=1}^p \frac{(x_j-EX_j)^2}{\sigma_j^2}] \ | $$
$$\approx  \lim_{|x_i| \rightarrow \infty} 3 \beta_{i,L_2}\frac{(x_i-EX_i)^2}{\sigma_i^2}=\infty;$$
last  approximation follows from (\ref{eq:sumIFderL2xrem}).  \hspace{3ex} $\Box.$

{\bf Proof of Proposition \ref{p:indexy}:} Follows from (\ref {eq:r2dify}) dividing both its sides by $h$ and taking the limit with $h$ converging to zero. $\hspace{2ex} \Box$

\begin{figure}
\resizebox{7in}{10in}{
\includegraphics{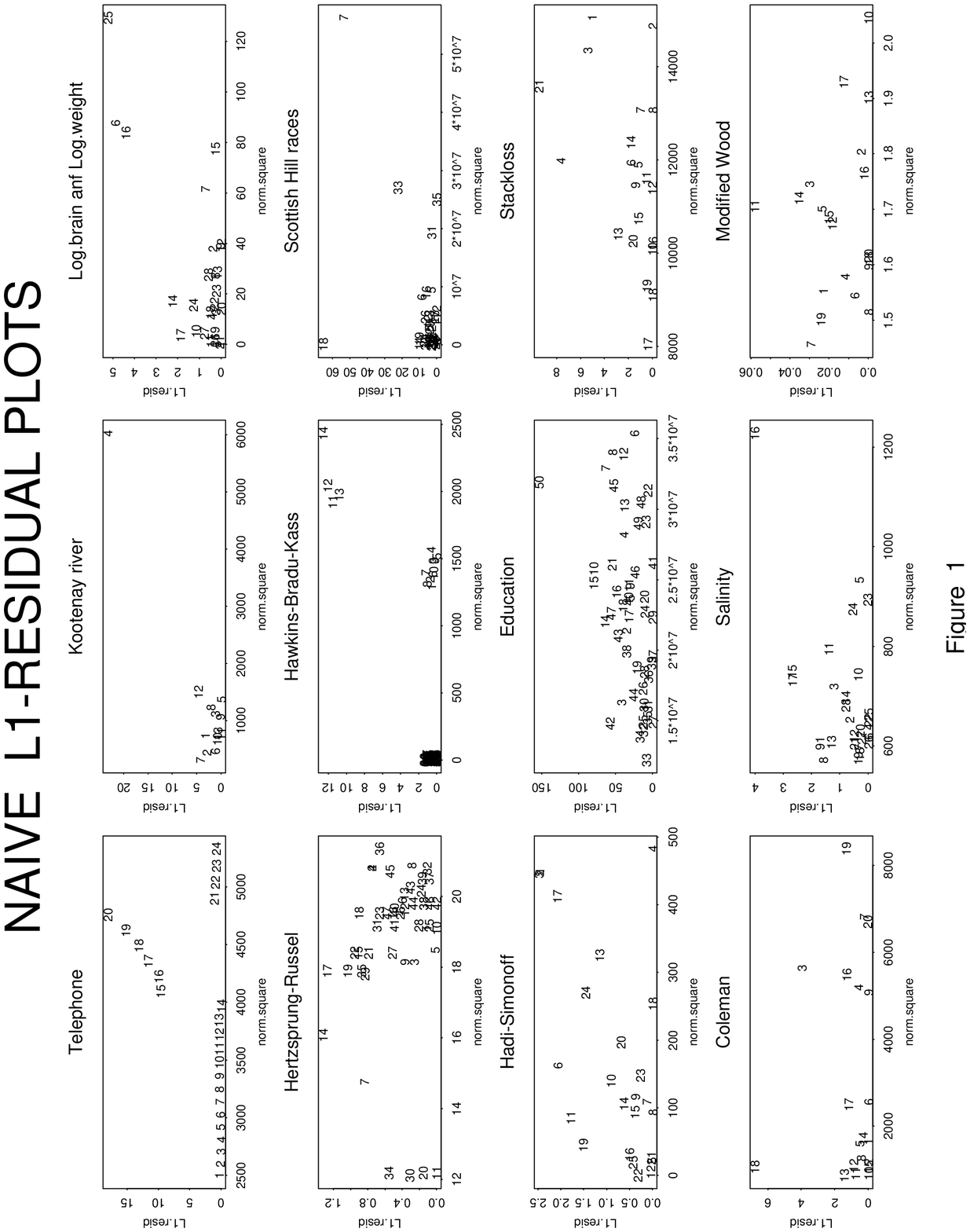}
}

\end{figure}

\begin{figure}

\resizebox{7in}{10in}{
\includegraphics{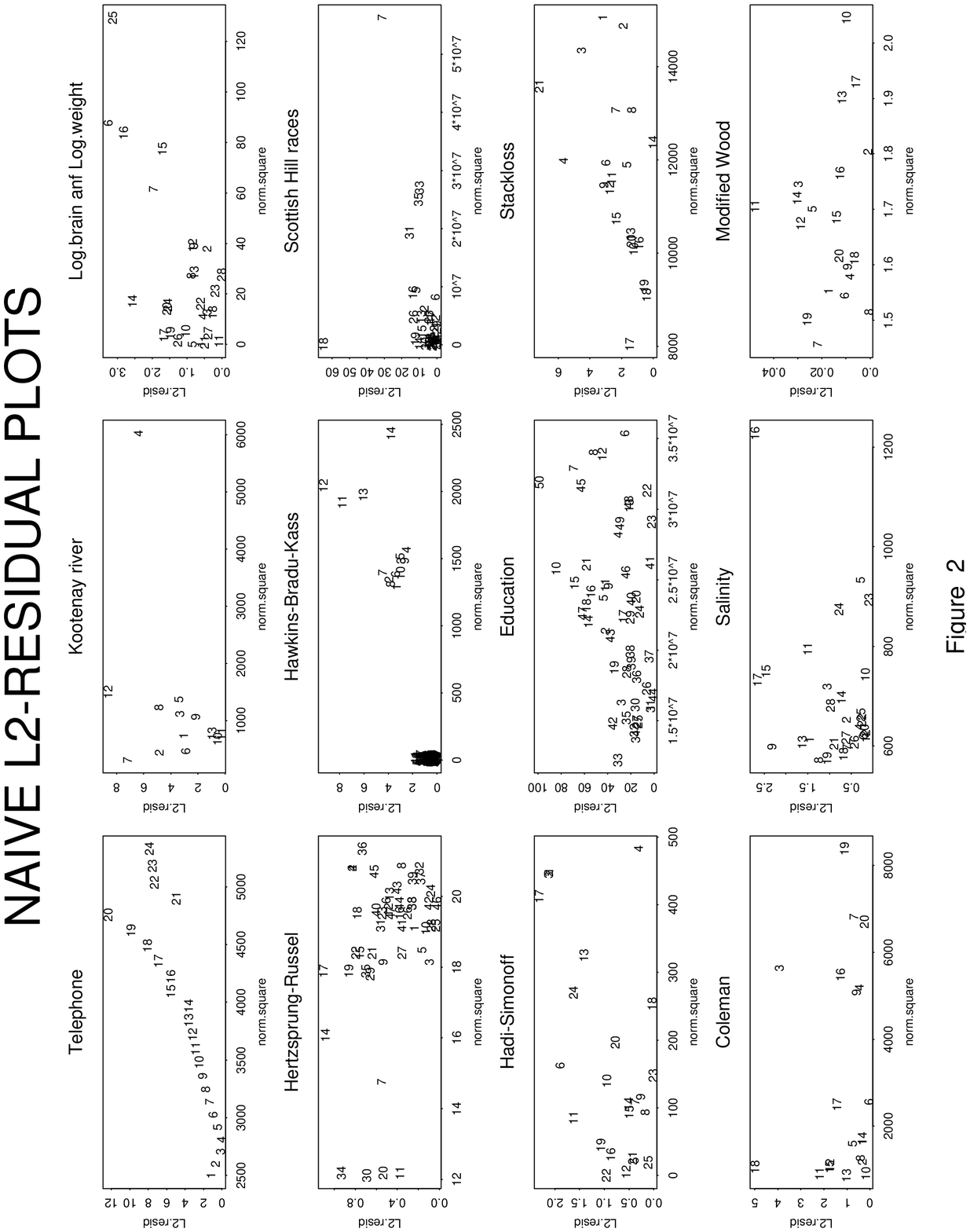}
}

\end{figure}

\begin{figure}

\resizebox{7in}{10in}{
\includegraphics{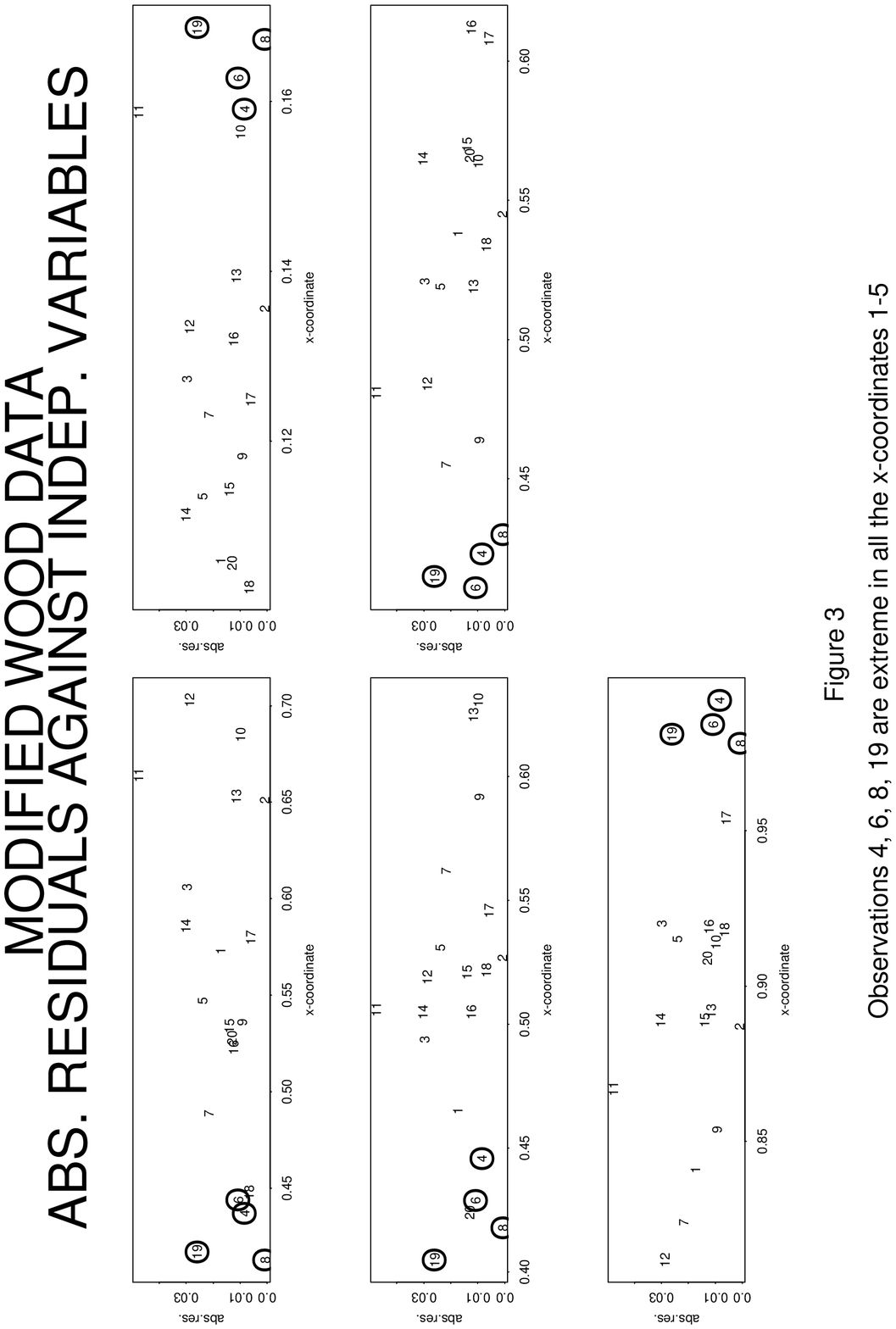}
}

\end{figure}

\end{document}